%% file: loc-unitary.tex
\documentclass[12pt]{amsart}
\usepackage{pb-diagram}
\usepackage{cmex12-2e}
\usepackage[6in]{amspapermacs}

\input loc-unimacs

\begin{document}

\title[Locally Inner Actions]{\boldmath
Locally Inner Actions on $\cox$-Algebras}

\author[Echterhoff]{Siegfried Echterhoff}
\address{Universit\"at-Gesamthochschule Paderborn \\
Fachbereich Mathematik-Informatik \\
Warburber Stra\ss e 100 \\
D-33095 Paderborn \\
Germany}
\email{echter@uni-paderborn.de}

\author[Williams]{Dana P. Williams}
\address{Department of Mathematics \\
Dartmouth College \\
Hanover, NH 03755-3551 \\
USA}
\email{dana.williams@dartmouth.edu}

\date{16 June 1997}

\ifdraft
\include{loc-uni-ab}\else
\input loc-uni-ab
\fi

\maketitle

\ifdraft

\include{loc-uni-intro}\else
\input loc-uni-intro
\fi
\ifdraft
\include{loc-uni-pre}\else
\input loc-uni-pre
\fi
\ifdraft
\include{loc-uni-lu}\else
\input loc-uni-lu
\fi
\ifdraft

\include{repr-gr}\else
\input repr-gr
\fi
\ifdraft
\include{brauer}\else
\input brauer
\fi
\ifdraft
\include{brauer1}\else
\input brauer1
\fi
\ifdraft
\include{ros}\else
\input ros
\fi

\includerefs

\end{document}

%% file: loc-unimacs.tex
\makeatletter
\@addtoreset{footnote}{section}
\@addtoreset{footnote}{subsection}
\makeatother

\theoremstyle{definition}
\newtheorem{remark}[thm]{Remark}
\newtheorem{ex}[thm]{Example}
\newtheorem{example}[thm]{Example}

\theoremstyle{definition}

\theoremstyle{plain}
\newtheorem*{thmnn}{Theorem}

\DeclareMathSymbol{\rtimes}{\mathbin}{AMSb}{"6F}

\def\modulefont#1{\mathsf{#1}}
\def\X{{\modulefont X}}

\def\ib{im\-prim\-i\-tiv\-ity bi\-mod\-u\-le}
\def\e{\xi}
\def\sec{\Gamma_0}
\def\s(#1){\sec(\e,#1)}
\def\hm{homomorphism}
\def\UM{UM}
\def\ai{\alpha^i}

\def\ui{u^i}
\def\uj{u^j}
\def\vi{v^i}
\def\vj{v^j}
\def\gij{\gamma_{ij}}
\def\gjk{\gamma_{jk}}
\def\gik{\gamma_{ik}}
\def\subab{\text{\normalfont ab}}
\def\gab{G_{\subab}}
\def\nab{N_{\subab}}
\def\hgab{\widehat{G}_{\subab}}
\def\hnab{\widehat{N}_{\subab}}
\def\sheaf#1{\mathcal {#1}}

\def\shgab{\widehat{\sheaf G}_{\subab}}
\def\CXtensor{\tensor_{C(X)}}
\def\Xtensor{\tensor_X}

\def\UB{U}
\def\H{\mathcal{H}}

\def\sme{\,\mathord{\mathop{\text{--}}\nolimits_{\relax}}\,}
\def\<{\langle}
\def\>{\rangle}
\let\ipscriptstyle=\scriptscriptstyle
\def\lipsqueeze{{\mskip -3.0mu}}
\def\ripsqueeze{{\mskip -3.0mu}}
\def\ipcomma{\mathrel{,}}
\newbox\ipstrutbox
\setbox\ipstrutbox=\hbox{\vrule height8.5pt
width 0pt}
\def\ipstrut{\copy\ipstrutbox}
\def\lip#1<#2,#3>{\mathopen{\relax_{\ipstrut\ipscriptstyle{
#1}}\lipsqueeze
\langle} #2\ipcomma #3 \rangle}
\def\blip#1<#2,#3>{\mathopen{\relax_{\ipstrut
\ipscriptstyle{ #1}}\lipsqueeze\bigl\langle} #2\ipcomma #3 \bigr\rangle}
\def\rip#1<#2,#3>{\langle #2\ipcomma #3
\rangle_{\ripsqueeze\ipstrut\ipscriptstyle{#1}}}
\def\brip#1<#2,#3>{\bigl\langle #2\ipcomma #3
\bigr\rangle_{\ripsqueeze\ipstrut\ipscriptstyle{#1}}}
\def\angsqueeze{\mskip -6mu}
\def\smangsqueeze{\mskip -3.7mu}
\def\trip#1<#2,#3>{\langle\smangsqueeze\langle #2\ipcomma #3
\rangle\smangsqueeze\rangle_{\ripsqueeze\ipstrut\ipscriptstyle{#1}}}
\def\btrip#1<#2,#3>{\bigl\langle\angsqueeze\bigl\langle #2\ipcomma
#3
\bigr\rangle
\angsqueeze\bigr\rangle_{\ripsqueeze\ipstrut\ipscriptstyle{#1}}}
\def\tlip#1<#2,#3>{\mathopen{\relax_{\ipstrut\ipscriptstyle{
#1}}\lipsqueeze \langle\smangsqueeze\langle} #2\ipcomma #3
\rangle\smangsqueeze\rangle}
\def\btlip#1<#2,#3>{\mathopen{\relax_{\ipstrut\ipscriptstyle{
#1}}\lipsqueeze
\bigl\langle\angsqueeze\bigl\langle} #2\ipcomma #3
\bigr\rangle\angsqueeze\bigr\rangle}
\def\LUproperty{LU}
\def\lu{\ifmmode\LUpropery\else$\LUproperty$\fi}

\def\fdbar{$\text{[FD]}\;\bar{}$}

\def\om{\omega}
\def\tg{\operatorname{tg}}
\def\RR{\mathbb{R}}

\def\TT{\mathbb{T}}
\def\ZZ{\mathbb{Z}}
\def\NN{\mathbb{N}}

\def\br{{\operatorname{Br}}}
\let\Br\br
\def\brg{\br_G}

\def\hA{\hat A}

\def\MyMath#1{\ifmmode #1\else$#1$\fi}
\def\prima{\MyMath{\Prim(A)}}
\def\cb{C^b}
\def\cbp{\cb\(\prima\)}

\def\M{M}
\def\UM{\mathcal{U}\M}
\def\ZM{\mathcal{Z}\M}
\def\ZUM{\mathcal{Z}\UM}
\def\cox{\MyMath{C_0(X)}}
\def\coxalg{\cox-algebra}
\DeclareMathOperator{\Glimm}{Glimm}
\def\gla{\ifmmode\Glimm(A)\else$\Glimm(A)$\fi}
\def\tcr{\tau_{\text{cr}}}
\def\tq{\tau_{\text{q}}}
\def\ewcr{\mathcal{CR}}
\def\CR{\MyMath{\ewcr}}
\def\U{U}
\def\LI{\MyMath{\mathcal{LI}}}
\def\cU{\boldsymbol{U}}
\def\tensor{\otimes}
\def\mtensor{\tensor_{\text{max}}}
\def\atmb{A\mtensor B}
\def\atxb{A\xtensor B}
\def\xtensor{\tensor_X}
\def\coxk{C_0(X,\K)}

\def\vec#1{\mathbf{#1}}
\def\vs{{\vec{s}}}
\def\E{\mathcal{E}}
\def\egx{\MyMath{\E_G(X)}}
\def\cxhtt{C\(X,H^2(G,\T)\)}
\def\igamma{\gamma^o}

%% file: loc-uni-ab.tex
%
%

\begin{abstract}
We make a detailed study of locally inner actions on \cs-algebras
whose primitive ideal spaces have locally compact Hausdorff complete
regularizations.  We suppose that $G$ has a representation group and
compactly generated abelianization $\gab$.  Then if the complete
regularization of $\prima$ is $X$, we show that the collection of
exterior equivalence classes of locally inner actions of $G$ on $A$ is
parameterized by the group $\E_G(X)$ of exterior equivalence classes of
\cox-actions of $G$ on $\coxk$.  Furthermore, we exhibit a group
isomorphism of
$\E_G(X)$ with the direct sum $H^1(X,\shgab)\oplus C\(X,H^2(G,\T)\)$.
As a consequence, we can compute the equivariant Brauer group
$\br_G(X)$ for $G$ acting trivially on $X$.
\end{abstract}

%% file: loc-uni-intro.tex
%
%

\section{Introduction}

One of the original motivations for the study of \cs-algebras arose
from the desire to understand the representation theory of locally
compact groups.  As is eloquently described in Rosenberg's survey
article \cite[\S3]{ros5}, the modern Mackey-Green machine shows that
to make further progress in this direction, it will be necessary to
have detailed knowledge of certain twisted crossed product
\cs-algebras --- either those studied by Green \cite{green1}, or more
generally, Busby-Smith twisted crossed products as studied by Packer
and Raeburn \cite{para1}.  In view of the Packer-Raeburn Stabilization
Trick \cite[Theorem~3.4]{para1}, it really suffices to consider
ordinary crossed products $A\rtimes_\alpha G$.  Of course any sort of
general classification of crossed products is out of the question;
however, as Rosenberg outlines in his ``Research Problem 1'' from
\cite{ros5}, it would be quite valuable and interesting to obtain
detailed information in  the
special case of an action with a ``single orbit type'' acting on a
continuous-trace \cs-algebra.    There is a considerable volume of
work in this direction --- for example, \cite{pr2,rw,rr,doir, echros,
ech4}  and other references cited in \cite[\S3]{ros5}.  Notice that
the case of nonvanishing Mackey obstructions is treated only in
\cite{echros, ech4},
and that in the majority of the published results it is assumed
that the group acting is abelian.

In this article, we consider a general family of dynamical systems
which include all spectrum fixing actions of a wide class of locally
compact groups acting on continuous-trace \cs-algebras.
Since the crossed product by any action of $G$ on a stable continuous-trace
\cs-algebra $A$ with single orbit type and  constant stabilizer $N$
can be decomposed, via the stabilization trick, into an spectrum fixing
action of $N$ and an action of $G/N$ on $A\rtimes_\alpha N$ with
$G/N$ acting freely on $\specnp{(A\rtimes_{\alpha} N)}$, a detailed
description of spectrum fixing actions and their crossed products provides a
major step towards a general solution to Rosenberg's ``Research
Problem 1.''  Following ideas developed in \cite{pr2, 90a, ckrw,
judymc, prw} (and others) we are going to describe our actions in
terms of topological invariants living in Moore's group cohomology and
certain sheaf cohomology groups.  In a forthcoming paper \cite{ew3}, we
will use these topological invariants to give a precise
bundle-theoretic description of the corresponding crossed products.
Since our methods require separability in a number of essential ways,
we will assume from the onset that {\em all automorphism groups are
second countable and that the \cs-algebras on which they act are
separable}.

Unlike much of the earlier work, we do not
assume that our groups are abelian, nor do we assume
that the associated Mackey obstructions vanish as in
\cite{pr2}, nor do we assume that they are constant as in
\cite{echros}. Moreover, it turns out that one
basic reason for actions on
continuous-trace \cs-algebras being more manageable than
actions on arbitrary \cs-algebras is that
for suitable $G$ (e.g., if the
abelianization $\gab=G/\overline{[G,G]}$ is compactly generated),
any spectrum fixing action  of $G$ on a continuous-trace \cs-algebra $A$ is
\emph{locally inner}.  (This follows from 
the proof of \cite[Corollary~2.2]{ros2}.)
This means that each point in the spectrum has
an open neighborhood $U$ such that the action on the corresponding ideal
$A_U$ of $A$ is {\em inner\/} in the sense that for all $s\in G$ there is a
unitary $u_s\in\UM(A_U)$ such that
$\alpha_s|_{A_U}=\Ad u_s$. Thus it is natural to try to classify locally
inner actions on arbitrary \cs-algebras rather than restricting ourselves
to actions on continuous-trace algebras. This turns out to
be possible for a large class of \cs-algebras, namely those
whose primitive ideal space $\Prim(A)$ has a second countable
locally compact complete regularization $X$ as described in
\secref{sec-1}.
This class of algebras includes all unital \cs-algebras, all
\cs-algebras with Hausdorff primitive ideal spaces, and all
quasi-standard \cs-algebras in the sense of Archbold and Somerset
\cite{as}.  If $X$ is given, then we will denote by $\CR(X)$ the class
of
\cs-algebras for which the complete regularization of $\Prim(A)$ is
homeomorphic to $X$.

From the point of view of dynamical systems, two actions $\alpha$ and
$\beta$ of $G$ on $A$ are considered to be ``the same'' if they are
\emph{exterior equivalent} (see \S 2 for the precise definition). For
instance, if $\alpha$ and $\beta$ are exterior equivalent, then the
corresponding crossed products are isomorphic.  An action $\alpha:G\to
\Aut(A)$ is exterior equivalent to the trivial action if and only if
it is \emph{unitary}; that is, if there is a strictly continuous
\emph{\hm} $u:G\to\UM(A)$ such that $\alpha_s=\Ad(u_s)$ for all $s\in
G$.  As will become apparent, the crucial step in describing locally
inner actions in general turns out to be describing the collection
$\egx$ of exterior equivalence classes of $\cox$-linear actions on
$A=\coxk$.  (Here and in the sequel, $\K$ will denote the compact
operators on a separable infinite dimensional Hilbert space $\H$.)  In
fact, $\egx$ is an abelian group.  It is a special case of
\cite[Theorem~3.6]{ckrw} that the collection of Morita equivalence
classes of \cox-systems $(A,G,\alpha)$, where $A$ is a
continuous-trace \cs-algebra with spectrum $X$, forms a group
$\brg(X)$ with respect to the operation
\begin{equation*}
[A,\alpha][B,\beta]:=[A\Xtensor B,\alpha\Xtensor\beta],
\end{equation*}
where $A\Xtensor B$ denotes the balanced tensor product $A\CXtensor B$
and $\alpha\Xtensor \beta$ is the action on the quotient induced by
$\alpha\tensor \beta$.
In this situation, two actions on $\coxk$ are Morita equivalent if and
only if they are exterior equivalent (\propref{prop-help}).
Therefore we can
identify \egx{} with the subgroup of $\brg(X)$ equal to the
kernel of the Forgetful \hm{}
$F:\brg(X)\to\br(X)\cong H^3(X;\Z)$, which maps the class of
$[A,\alpha]$ to the Dixmier-Douady class $\delta(A)$ of $A$ in the
third cohomology of $X$.  Since $G$ acts trivially here, $F$ is
clearly surjective and admits a natural splitting map so that the
equivariant Brauer group
$\brg(X)$ is isomorphic to $ \egx\oplus H^3(X;\Z)$.

To fix ideas and notation, consider the
(well known) case $X=\{pt\}$.
If $\beta:G\to \Aut\(\K(\H)\)$ is an action, then, from the
short exact sequence
of polish groups
\[1\arrow{e} \TT\arrow{e} \U(\H)\arrow{e,t}{\Ad} \Aut(\K)\arrow{e} 1\]
we obtain an obstruction $[\om_{\beta}]\in
H^2(G,\TT)$ (Moore's group cohomology) to lifting
$\beta$ to a strongly continuous homomorphism $V:G\to \U(\H)$;
thus, $[\om_{\beta}]=0$ if and only if $\beta$ is unitary.
This obstruction is
called the {\em Mackey obstruction\/} of $\beta$.
Notice that if $\om\in Z^2(G,\TT)$, then $[\om]=[\om_{\beta}]$
if and only if there exists an $\om$-representation
$V:G\to\U(\H)$ which implements $\beta$; that is,
$V$ is a Borel map satisfying
$V_e=1$, $V_sV_t=\om(s,t)V_{st}$ for $s,t\in G$, and
$\beta_s=\Ad V_s$ for all $s\in G$.
Using this  it is not hard to show that
$[\beta]\mapsto[\om_{\beta}]$ is indeed a group isomorphism
between $\E_G(\{pt\})\to H^2(G,\TT)$; surjectivity follows
from the fact that for each $\om\in Z^2(G,\TT)$ there
is at least one $\om$-representation, namely the left
regular $\om$-representation on $L^2(G)$ defined by
$(L^\om_s\xi)(t)=\om(s,s^{-1}t)\xi(s^{-1}t)$.

Let $\beta:G\to \Aut\(C_0(X,\K)\)$ be a $C_0(X)$-linear action
and let $[\om_x]\in H^2(G,\TT)$ denote the Mackey-obstruction
for the induced automorphism group $\beta^x$ acting on the fibre over
$x$. The map $\varphi^{\beta}:X\to H^2(G,\TT)$ given by $
\varphi^{\beta}(x)=[\om_x]$ is continuous (when
$G$ is compactly generated, this follows from \cite[Lemma 3.3]{doir}; for
the general case see \lemref{lem-fix} below.)
Since the evaluation map $\E_G(X)\to\mathcal
E_G(\{x\})$ sending $ [\beta]\mapsto [\beta^x]$ is a homomorphism for
each $x\in X$, it follows from  the discussion above
that $[\beta]\mapsto \varphi^{\beta}$ is a group
homomorphism
$$\Phi:\E_G(X)\to C\(X,H^2(G,\TT)\).$$
Notice that $\varphi^{\beta}=0$ if and
only if
each irreducible representation of $C_0(X,\K)$ can be extended
to a covariant representation of
$\big(C_0(X,\K),G,\beta\big)$; that is, if and only if $\beta$ is
{\em pointwise unitary}. Thus $\ker \Phi$
 consists of all exterior equivalence
classes of pointwise unitary actions of $G$ on $C_0(X,\K)$.

In \cite{judymc} (see also \cite{prw}), Packer observed that if $G$ is an
elementary abelian group (i.e., $G$ is of the form
$\R^n\times\T^m\times\Z^k\times F$,
$F$ finite abelian), then
$\Phi$ is surjective and admits a splitting map.  Moreover, a theorem of
Rosenberg (see \thmref{thm-ros}) implies that under these hypotheses,
$\ker \Phi$ coincides with the (equivalence classes of) locally
unitary actions.  Therefore the Phillips-Raeburn obstruction map (see
\corref{cor-point}) gives an isomorphism of $\ker \Phi$ with the
isomorphism classes of principal $\widehat G$-bundles over $X$, or
equivalently, with the sheaf cohomology group $H^1(X,\widehat{\sheaf
G})$.  (If $G$ is an abelian group, we will use the corresponding
caligraphic letter $\sheaf G$ to denote the sheaf of germs of
continuous $G$-valued functions.)  Thus as abelian groups, $\egx\cong
H^1(X,\widehat{\sheaf G}) \oplus \cxhtt$, and $\brg(X)\cong
H^1(X,\widehat{\sheaf G}) \oplus \cxhtt \oplus H^3(X;\Z)$.  Notice
that the critical steps in Packer's argument are to find a splitting
map for $\Phi$ and to identify $\ker \Phi$ using Rosenberg's theorem.

Our first main result is to produce a splitting map for $\Phi$ in the
case that $G$ has a representation group in the sense of Moore (see
Definition~\ref{def-3.1}).  Such groups are called \emph{smooth} and
comprise a large class of locally compact groups including all compact
groups, all discrete groups, and all compactly generated abelian
groups (see Remark~\ref{rem-rep} and \corref{corcompgen}).  Smooth
groups $G$ have the property that $H^2(G,\T)$ is locally compact and
Hausdorff.  Thus if $\gab$ is compactly generated, then $G$ satisfies
the hypotheses of Rosenberg's theorem and allows us to identify
$\ker \Phi$ with locally unitary actions of $G$.  A suitable
modification of the Phillips-Raeburn theory (see \secref{sec-2}) gives
an isomorphism of $\ker \Phi$ with $H^1(X,\shgab)$
(\corref{cor-point}).   Thus we obtain the following result.

\begin{thmnn}
[{\thmref{brauer} \& \corref{cor-brauer}}]
Suppose that $G$ is smooth and that $\gab$ is compactly generated.
Then
\begin{equation*}
\E_G(X)\cong H^1(X,\shgab)\oplus \cxhtt,
\end{equation*}
and for any trivial $G$-space $X$,
\begin{equation*}
\brg(X)\cong H^1(X,\shgab)\oplus \cxhtt \oplus H^3(X;\Z).
\end{equation*}
\end{thmnn}
This gives
new information even if $G$ is abelian, since
by Corollary~\ref{corcompgen} our result applies not only to
elementary abelian groups, but to all
second countable compactly generated abelian groups.

Our final result is a consequence of the above theorem and the
observation that if $\alpha$ is a locally inner action of a smooth
group $G$ on $A\in\CR(X)$, then there is a $[\gamma]\in\egx$ such that
$\alpha\Xtensor \gamma$ is locally unitary on $A\Xtensor\coxk$
(\propref{prop-locuni}).

\begin{thmnn}
[{\thmref{thm-locinner}}]
Suppose that $A\in \CR(X)$ and $\alpha:G\to\Aut(A)$ is a locally inner
action of a smooth group on $A$.  If $\gab$ is compactly
generated, then there is a unique $[\beta^\alpha]\in\egx$ such that
$\alpha\Xtensor
\id$ is exterior equivalent to $\id\Xtensor \beta^\alpha$ on
$A\Xtensor \coxk$.  In fact, the map $[\alpha]\mapsto [\beta^\alpha]$
is a well-defined injective map from the collection of exterior
equivalence classes of locally inner actions on $A$ to \egx.  This
correspondence is bijective if $A$ is stable.
\end{thmnn}

It follows that if $A\in\CR(X)$ and if $(A,G,\alpha)$ is locally
inner, then the crossed product $A\rtimes_\alpha G$ is Morita
equivalent to one of the special form $\big(A\Xtensor
\coxk\big)\rtimes_{\id\Xtensor\beta} G$, where $\beta$ is in \egx.
Having this, it is possible to describe the crossed product in
terms of the invariants associated to $\beta$ and a representation
group for $G$; this we do in \cite{ew3}.

Our work is organized as follows.  \secref{sec-1} is devoted to some
preliminary results on \cox-algebras and on the complete regularization
of \prima.  In \secref{sec-2} we consider locally unitary actions on
$\CR(X)$-algebras, and extend the Phillips-Raeburn
classification scheme to this setting.  Since smooth groups play such
an important r\^ole in the sequel, we devote \secref{sec-3} to
developing some basic results about
representation groups.  The chief result connecting representation
groups to the splitting of $\Phi$ is the characterization of smooth
groups given in \lemref{cocycle}.  In \secref{sec-4}, we prove the
first of our main results (\thmref{brauer}) which describes \egx.  Our
description of locally inner actions is given in \secref{sec-5}.

Since Rosenberg's theorem plays a key r\^ole, we provide a discussion
of possible extensions of his theorem in \secref{sec-appendix}.  We
give examples which show that the hypotheses are sharp --- that is,
the major assumptions that $H^2(G,\T)$ is Hausdorff and that $\gab$ is
compactly
generated are both necessary in general.  On the other hand,
we also show that
Rosenberg's theorem holds for a strictly larger class of groups if we
restrict ourselves to actions on continuous-trace \cs-algebras with
locally connected spectrum (\thmref{append}).  This class of groups
contains all connected nilpotent Lie groups and all
\fdbar{} groups (\corref{cor-FD}) --- a
class of groups which contains all known examples of groups $G$ for which
$\cs(G)$ is a continuous-trace
\cs-algebra.

%% file: loc-uni-pre.tex
%
%

\section{Preliminaries}\label{sec-1}

If $A$ is a \cs-algebra, then we will write $\prima$ for the space
of primitive ideals of $A$ with the Jacobson topology.  This topology
is badly behaved in general and may satisfy only the $T_0$-axiom for
separability.  On the other hand, \prima{} is always locally
compact\footnote{We do not require that compact or locally compact
spaces be Hausdorff.}, and \prima{} is second countable whenever $A$
is separable \cite[\S3.3]{dix}.
The Jacobson topology on \prima{} not
only describes the ideal structure of $A$, but also allows us to
completely describe the center $\ZM(A)$ of the multiplier algebra
$\M(A)$ of $A$.  If $a\in A$, then we will write
$a(P)$ for the image of $a$ in the quotient $A/P$, then
the Dauns-Hofmann Theorem allows us to identify
$\cbp$ with $\ZM(A)$ as follows: if $f\in\cbp$ and if $a\in A$, then
$f\cdot a$ is the unique element of $A$ satisfying $(f\cdot a)(P) =
f(P)a(P)$ for all $P\in\prima$,
and every element of $\ZM(A)$ is of this form  (cf.,
\cite[Corollary~4.4.8]{ped} or \cite{may}).  Note that
$A$ is a nondegenerate central Banach $\cbp$-module.

Since the topology on \prima{} can be awkward to deal with, a natural
alternative is to use the following definition.

\begin{definition}
Suppose that $X$ is a locally compact Hausdorff space.  A
\emph{\coxalg} is a \cs-algebra $A$ together with a $*$-\hm{}
$\Phi_A:\cox\to\ZM(A)$ which is \emph{nondegenerate} in the
sense that
\[
\Phi_A\(\cox\)\cdot A:=\sp\set{\Phi_A(f)a:\text{$f\in \cox$ and $a\in A$}}
\]
is dense in $A$.
\end{definition}

\coxalg s have enjoyed a considerable amount of attention recently and
there are a number of good treatments available
\cite{blanchard,blanchard2,may1}.  We recall some of the basic
properties here.

If $(A,\Phi_A)$ is a \coxalg, then there is a continuous map
$\sigma_A:\prima\to X$ such that $\Phi_A(f)=f\circ\sigma_A$.  (Here
and in the sequel, we \emph{identify} $\ZM(A)$ with $\cbp$ via the
Dauns-Hofmann Theorem.)
As the converse is clear, $A$ is a \coxalg{} if and only if there is a
continuous map from $\prima$ to $X$.  We will usually suppress
$\Phi_A$ and $\sigma_A$ and write $f\cdot a$ in place of $\Phi_A(f)a$
or $(f\circ\sigma_A)\cdot a$.  Notice that $A$ is a nondegenerate
central Banach \cox-module satisfying
\begin{equation}
\label{eq-//}
(f\cdot a)^*=a^*\cdot \bar f.
\end{equation}
Furthermore, any nondegenerate central \cox-module satisfying
\eqref{eq-//} is a \coxalg.

Suppose that $U$ is open in $X$ and that $J$ is the ideal of functions
in \cox{} vanishing off $U$.  Then the Cohen factorization theorem
(\cite{cohen}, \cite[Proposition~1.8]{blanchard}) implies that
\begin{equation*}
\overline{J\cdot A}:=\overline\sp\set{f\cdot a:\text{$f\in J$ and
$a\in A$}}
=\set{f\cdot a:\text{$f\in J$ and
$a\in A$} }.
\end{equation*}
(The point being that it is unnecessary to take either the closure or the
span in the final set.)
Anyway, $J\cdot A$ is an ideal of $A$ which will be denoted by $A_U$.
For each $x\in X$, we write $A(x)$ for the quotient of $A$ by
$A_{X\setminus\set x}$.  If $a\in A$, then we write $a(x)$ for the
image of $a$ in $A(x)$.  We refer to $A(x)$ as \emph{the fibre of $A$
over $x$}.  Notice that it is possible that $A(x)=\set0$.
Even so, we often view elements of $A$ as ``fields'' in
$\bigoplus_{x\in X} A(x)$.  This point of view is justified by the
following.

\begin{lem}
[{\cite{blanchard,may1}}]
Suppose that $A$ is a \coxalg.  For each $a\in A$, the map $x\mapsto
\|a(x)\|$ is upper semicontinuous; that is, $\set{x\in
X:\|a(x)\|\ge\epsilon}$ is closed for all $\epsilon\ge0$.
Furthermore,
\[
\|a\|=\sup_{x\in X}\|a(x)\|.
\]
\end{lem}

\begin{remark}
The map
$x\mapsto
\|a(x)\|$ is continuous for all $a\in A$ if and only if
$\sigma_A$ is
open \cite{lee2,may1}.  In this case,
$A$ is the section algebra of a
\cs-bundle over the image of $\sigma_A$ \cite[\S1]{fell4}.
\end{remark}

\begin{remark}
Notice that if $A$ is a \coxalg, then each $m\in\M(A)$ defines a
multiplier $m(x)\in \M\(A(x)\)$.  If $m\in\M(A)$
and $a\in A$, then $ma(x)=m(x)a(x)$.
\end{remark}

Given a \cs-algebra $A$, it is natural to look for a nice space $X$
which makes $A$ a \coxalg.   Since $X$ will be the image of \prima{}
by a continuous map, it is reasonable to look for a
``Hausdorffication
of $\prima$''.  Regrettably, there are a
horrifying number of alternatives to chose from (cf., e.g.,
\cite[Chap.~III
\S3]{dauns-hofmann}).  For our purposes, the appropriate
notion is the \emph{complete regularization}.  If $P$ and $Q$
belong to \prima, then we define
$P\sim Q$ if $f(P)=f(Q)$ for all $f\in \cbp$.  Then $\sim$ is an
equivalence relation and the set $\prima/\!\!\sim$ is denoted by
$\gla$
\cite{as}.  If we give \gla{} the weak topology $\tcr$ induced by the
functions in $\cbp$, then $\(\gla,\tcr\)$ is a completely regular
space \cite[Theorem~3.7]{gill-jer}.   The quotient map
$q:\prima\to\gla$ is called the \emph{complete regularization map}.
It is not clear that $\tcr$ coincides with the quotient
topology\footnote{These topologies do differ in general
\cite[3J.3]{gill-jer}; however, we know of no examples where they
differ for $\gla$.} $\tq$
on $\gla$, although one certainly has $\tcr\subseteq \tq$.  In
particular, $q$ is continuous; moreover the map $f\mapsto f\circ q$ is
an isomorphism of $\cb\(\gla)$ and $\cbp$
\cite[Theorem~3.9]{gill-jer}.
Furthermore, $\tcr$ is the only completely regular topology on \gla{}
such that the functions induced by $\cbp$ are continuous
\cite[Theorem~3.6]{gill-jer}.

Here it will be necessary to have the complete regularization $\(\gla,\tcr\)$
be locally compact.  Regrettably, this can
fail to be the case
\cite[Example~9.2]{dauns-hofmann}.   Even if the complete
regularization is locally compact, we have been unable to show that it
must be second countable if $A$ is separable.  Consequently, we must
include both these assumptions in our applications.

\begin{definition}
We will call a separable \cs-algebra a \emph{\CR-algebra} if the complete
regularization $X:=\(\gla,\tcr\)$ of \prima{} is a second countable
locally compact Hausdorff space.  If $X$ is a second countable
locally compact Hausdorff space, then we will write $\CR(X)$ for the
collection of \CR-algebras with complete regularization homeomorphic
to $X$.
\end{definition}

Despite the pathologies mentioned above, the class of \CR-algebras is
quite large.  It clearly contains all separable $C^*$-algebras with
Hausdorff primitive ideal space \prima{}. If $A$ is unital,
then \prima{} is compact.  Since the complete regularization map is continuous,
\gla{} is compact.  Since
$\cbp=\cb\(\gla\)$ is actually a closed subalgebra of $A$ in this
case, $\cb\(\gla\)$ is separable and $\gla$ is second
countable\footnote{There is an embedding of \gla{} into
$\specnp{\cb\(\gla\))}$ which is the Stone-\v Cech compactification
$\beta\(\gla\)$.
Since subset of a
locally compact Hausdorff space is locally compact if and only if it
is open in its closure,
$\gla$ is locally compact exactly when it is open in
its Stone-\v Cech compactification.}.
Thus
\emph{every unital \cs-algebra is a \CR-algebra}.
Another large class of \CR-algebras is provided by the
\emph{quasi-standard \cs-algebras} studied in \cite{as}.
Recall that a \cs-algebra is called quasi-standard if (1)~defining
$P\approx Q$ when $P$ and $Q$ cannot be separated by open sets in
\prima{} is an equivalence relation on \prima, and~(2) the
corresponding quotient map is open.  If $A$ is quasi-standard, then
$\sim$ and $\approx$ coincide and $A$ is \CR{}
\cite[Proposition~3.2]{as}.  In fact a number of interesting group
\cs-algebras
turn out to be quasi-standard \cite{arch-kan,kan-sch-tay}.

%
%

Let $\M(A)$ be the multiplier algebra of $A$.
Recall that a net $\set{T_i}$ converges to $T$ in the strict topology
on $\M(A)$ if
and only if $T_ia\to Ta$ and $T_i^*a\to T^*a$ for all $a\in A$.  If
the net is bounded, then it suffices to take $a$ in the unit ball of
$A$.  In fact, if each $T_i$ is unitary, then it suffices to check only
that $T_ia\to Ta$ for $a\in A$ with $\|a\|\le 1$.
Consequently if $A$ is separable, then the unitary group
$\UM(A)$ is a second countable topological group in the strict
topology
which admits a
complete metric (compatible with the topology).  That is, $\UM(A)$ is
a Polish group.  Since $\ZUM(A)$ is closed in $\UM(A)$, it too is a
Polish group.

For notational convenience, let $X=\(\gla,\tcr)$ be the complete
regularization of $\prima$.  Then we can identify $\ZM(A)$ with
$\cb\(X)$, and $\ZUM(A)$ with $C(X,\T)$.  However it is not
immediately obvious how to describe the strict topology on $C(X,\T)$.
Our next result says that when $X$ is second countable and locally
compact, then
the strict topology on $C(X,\T)$ coincides with the compact-open
topology (the topology of uniform convergence on compacta).

\begin{lem}
\label{lem-co}
Suppose that $X$ is a second countable locally compact Hausdorff space
and that $A\in\CR(X)$.  Then $\ZUM(A)$ with the strict topology is
homeomorphic to $C(X,\T)$ with the compact-open topology.
\end{lem}
\begin{remark}
The lemma holds for $X=(\gla,\tcr)$ whenever $C(X,\T)$ is a Polish
group in the compact-open topology.
In general,
$X$ is a
$\sigma$-compact, completely regular space.  If $\tq=\tcr$, then $X$
is compactly generated (or a $k$-space) by \cite[43H(3)]{willard}, and
at least in a compactly generated space, the limit of continuous
functions in the compact-open topology is continuous.  In order that
$C(X,\T)$ be metric, it seems to be necessary that $X$ be
``hemicompact'' \cite[43G(3)]{willard}.  In any case, if $X$ is
hemicompact, then the compact open topology is metric and complete.
In this case, $C(X,\T)$ is Polish, at least
\emph{provided} that
$X$ is second countable --- so that $C(X,\T)$ is separable.
However we have been unable to show that $X$ is second countable ---
even if $X$ is locally compact.
\end{remark}
\begin{proof}
Suppose that $f_n\to f$ uniformly on compacta in $C(X,\T)$.  If $a\in
A$ is nonzero and if $\epsilon>0$, then the image of
\[
C=\set{P\in\prima:\|a(P)\|\ge\epsilon/2}
\]
is compact in $X$ \cite[\S3.3]{dix}.  Thus there is an $N$ such that
$n\ge N$ implies that $|f_n\(q(P)\)-f\(q(P)\)|<\epsilon/\|a\|$ for all
$P\in C$.  If $P\notin C$, then for all $n$,
\[
\|f_n\cdot a(P)-f\cdot
a(P)\|\le \epsilon.
\]
  This proves that $f_n\to f$ strictly.

Since $X$ is second countable and locally compact, there is a sequence
of compact sets $\set{K_n}$ in $X$ such that $X=\bigcup_n K_n$ and
such that every compact set $K$ in
$X$ is contained in some $K_n$.  Then if $\set{V_n}$ is a countable
basis for the topology of $\T$, we get a sub-basis $\set{U_{n,m}}$
for the
compact-open topology on $C(X,\T)$ by setting
\[
U_{n,m}:=\set{f\in C(X,\T):f(K_n)\subseteq V_m}.
\]
It follows that $C(X,\T)$ is second countable in the compact-open
topology.  Using the $K_n$'s, it is easy to construct a complete
metric on $C(X,\T)$ compatible with the compact-open topology.  Thus,
$C(X,\T)$ is a Polish group in the compact open topology.  Since the
first part of the proof shows that the identity map is continuous from
the compact-open topology to the strict topology, the result follows
from the Open Mapping Theorem \cite[Proposition~5(b)]{moore3}.
\end{proof}

An automorphism $\alpha$ of a \cs-algebra $A$
is called \emph{inner} if there is a $u\in\UM(A)$ such that $\alpha =
\Ad(u)$.  (Recall that $\Ad(u)(a):=uau^*$.)  An action
$\alpha:G\to\Aut(A)$ is called \emph{inner} if $\alpha_s$ is inner for
each $s\in G$.  An action is called \emph{unitary} if there is a
strictly continuous \hm{} $u:g\to\UM(A)$ such that $\alpha_s=\Ad(u_s)$
for all $s\in G$.
Unitary actions are considered trivial; for example, if $\alpha$ is
unitary, then the crossed product $A\rtimes_\alpha G$ is isomorphic to
$A\mtensor \cs(G)$.  Also two actions $\alpha:G\to\Aut(A)$ and
$\beta:G\to\Aut(A)$ are called \emph{exterior equivalent} if there is
a strictly continuous map $w:G\to\UM(A)$ such that
\begin{enumerate}
\item
$\alpha_s(a) =
w_s\beta_s(a) w_s^*$ for all $a\in A$ and $s\in G$, and
\item
for all
$s,t\in G$, $w_{st}=w_s\beta_s(w_t)$.
\end{enumerate}
In this event, we call $w$ a $1$-cocycle.  Actions
$\alpha:G\to\Aut(A)$ and $\beta:G\to\Aut(B)$ are called \emph{outer
conjugate} if there is a $*$-isomorphism $\Phi:A\to B$ such that
$\beta$ and $\Phi\circ\alpha\circ\Phi^{-1}$ are exterior equivalent.

Although unitary actions are trivial from the point of view of
dynamical systems, 
inner actions can be quite interesting.  Another class of interesting
actions are those which are locally inner or even unitary.

\begin{definition}
\label{def-2.3}
Let $X$ be a second countable locally compact Hausdorff space and $G$
a second countable locally compact group.
Suppose that $A\in\CR(X)$ and that $\alpha:G\to\Aut(A)$ is an action.
Then $\alpha$ is called \emph{locally unitary} (\emph{locally
inner})
if every point in $X$
has a \nbhd{} $U$ such that $A_U$ is invariant under $\alpha$ and the
restriction $\alpha^U$ of $\alpha$ to $A_U$ is unitary (inner).
\end{definition}

\begin{remark}
If $A$ has Hausdorff spectrum $X$, then the
above definition coincides with the usual notion of a locally unitary
action (cf., \cite[\S1]{pr2} and \cite[\S1]{ros2}).
\end{remark}

Recall from \lemref{lem-co} that if $A\in \CR(X)$ then the group
$\ZUM(A)$ of $A$, equipped with the strict
topology, is isomorphic to the polish group $C(X,\TT)$
equipped with the compact-open topology. Thus we obtain a short exact
sequence
of polish groups
$$1\arrow{e} C(X,\TT)\arrow{e} \UM(A)\arrow{e} \Inn(A)\arrow{e} 1.$$
If $\alpha:G\to \Aut(A)$ is inner, then $\alpha$ defines a continuous
\hm{} of $G$ into $\Inn(A)$ with its Polish topology
\cite[Corollary~0.2]{rr}.
Thus, we can choose a
Borel map $V:G\to \UM(A)$ such that $V_e=1$ and such that
$\alpha=\Ad V$. (For example, if $c:\Inn(A)\to \UM(A)$ is a Borel
section such that $c(\id)=1$, then $V=c\circ \alpha$ will do
the job.) Then $V$ determines a Borel cocycle $\sigma\in
Z^2\(G,C(X,\TT)\)$ via the equation $V_sV_t=\sigma(s,t)V_{st}$
for $s,t\in G$. The class $[\sigma]\in H^2\(G, C(X,\TT)\)$
only depends on $\alpha$ and is the unique obstruction for
$\alpha$ being unitary (see \cite[Corollary 0.12]{rr}). In what
follows we  will refer to $V$ as a \emph{$\sigma$-homomorphism} of
$G$ into
$\UM(A)$ which implements~$\alpha$.

\begin{remark}
\label{rem-fix}
Suppose that $\alpha:G\to\Aut\(\coxk\)$ is an inner automorphism group
which is implemented by a $\sigma$-\hm{} as above.  Then the Mackey
obstruction for the induced action $\alpha^x$ on the fibre over $x$ is
the class of $\sigma(x)$, where $\sigma(x)$ is the cocycle in
$Z^2(G,\T)$ obtained by evaluation at $x$:
$\sigma(x)(s,t):=\sigma(s,t)(x)$.
\end{remark}

%% file: loc-uni-lu.tex
%
%

\section{Locally Unitary Actions}\label{sec-2}

In this section we want to see that the Phillips-Raeburn
classification scheme for locally unitary actions
of abelian groups can be extended to
all second countable locally compact groups acting on
\CR-algebras.

Suppose that $A\in\CR(X)$ and that $\alpha:G\to\Aut(A)$ is locally
unitary.  Then there is a cover $\cU=\set{U_i}_{i\in I}$ such that
$\alpha^{U_i}= \Ad(u^i)$ for strictly continuous \hm s
$u^i:G\to\UM(A_{U_i})$.  For convenience, we will write $U_{ij}$ for
$U_i\cap U_j$, $A_{ij}$ in place of $A_{U_{ij}}$, $\alpha^{ij}$ for
$\alpha^{U_{ij}}$, and $u^{ij}$ for $(u^i)^{U_{ij}}$.  Even though
$u^{ij}\not= u^{ji}$, both $u^{ij}$ and $u^{ji}$ implement
$\alpha^{ij}$.  It follows that for each $s\in G$,
$u_s^{ij}(u_s^{ji})^*$ belongs to $\ZUM(A_{ij})$.  In order to
identify $\ZUM(A_{ij})$ with $C(U_{ij},\T)$ (with the compact-open
topology), we need the following lemma.

\begin{lem}
\label{lem-crrestr}
Suppose that $A\in\CR(X)$ and that $U$ is open in $X$.  Then $A_U \in
\CR(U)$.
\end{lem}

\begin{proof}
Let $q:\prima\to X$ be the quotient map.  Then we can identify
$\Prim(A_U)$ with $q^{-1}(U)$, and $\Glimm(A_U)$ is the
quotient of the latter with topology induced by $C^b\(q^{-1}(U)\)$.
We have to show that $\Glimm(A_U)$ can be identified with $U$ with the
relative topology.  However since any $f\in\cbp$ restricts to an
element of $\cb\(q^{-1}(U)\)$, it is clear that $P\sim Q$ in
$q^{-1}(U)$ implies that $P\sim Q$ in \prima.  On the other hand,
suppose that $P\sim Q$ in \prima{} and that $f\in \cb\(q^{-1}(U)\)$.
Since $X$ is locally compact and $U$ is an open \nbhd{} of
$q(P)=q(Q)$, there is a $g\in C_c^+(X)$ with $g\(q(P)\)=1$ and
$\supp(g) \subseteq U$.  Therefore we may view $h= f(g\circ q)$ as an
element of $\cbp$.  Since $h(P)=h(Q)$ by assumption, we must have
$f(P)=f(Q)$.  Thus the two equivalence relations coincide on
$q^{-1}(U)$ and we can identify $U$ with $\Glimm(A_U)$ at least as a
set.

Let $\tau_r$ be the relative topology on $U$.  A similar
argument to that in the previous paragraph shows that any element of
$\cb\(q^{-1}(U)\)$ agrees at least locally with an element of $\cbp$.
Thus $\cb(U,\tau_r)$ and $\cb(U,\tcr)$ coincide.  Since both
topologies on $U$ are completely regular (Hausdorff) topologies, and
therefore are determined by the zero sets of $\cb(U)$
\cite[Theorem~3.7]{gill-jer}, the topologies coincide.
\end{proof}

Now the previous lemma allows us to conclude that
$A_{ij}\in\CR(U_{ij})$ so that we may identify $\ZUM(A_{ij})$ with
$C(U_{ij}, \T)$ as claimed.  Notice that $u^{ij}(x)=u^i(x)$.
Since
$C(U_{ij},\T)$ has the compact open
topology and $s\mapsto u_s^{ij}(u_s^{ji})^*$ is continuous, it follows
that
\[
(x,s)\mapsto u_s^{i}(x)u_s^{j}(x)^*
\]
is jointly continuous from $U_{ij}\times G$ to $\T$.
In particular, for each $x\in U_{ij}$, $s\mapsto
u_s^{i}(x)u_s^{j}(x)^*$ is a continuous character $\gamma_{ij}(x)$
on $G$, and
\begin{equation}
\label{eq-*}
u_s^{i}(x)=\gij(x)(s)u^{j}_s(x).
\end{equation}
Since any character has to kill the closure of the  commutator subgroup
$[G,G]$, we will always view $\gamma_{ij}(x)$ as a character on the
abelian group $\gab:=G/\overline{[G,G]}$.
The group $\gab$ is a locally compact abelian group usually called the
\emph{abelianization} of $G$.  Notice that the joint continuity
implies that the functions $\gamma_{ij}:U_{ij}\to\hgab$ are continuous
when $\hgab$ is given the usual locally compact dual topology (of
uniform convergence on compacta).  A straightforward computation using
the definition of the $\gamma_{ij}$'s shows that if $x\in U_{ijk}$,
then
\[
\gij(x)\gjk(x)=\gik(x).
\]
Thus the collection $\gamma=\set{\gij}$ defines a $1$-cocycle
in\footnote{If $G$ is an abelian group,
we use the caligraphic letter $\sheaf G$ to denote the corresponding
sheaf of germs of continuous $G$-valued functions on $X$.}
$Z^1(\cU,\shgab)$ and therefore a class $\zeta$ in $H^1(X,\shgab)$.
We claim this class depends only on
$(A,G,\alpha)$.  Suppose we had taken a different cover
$\set{V_j}_{j\in J}$ and \hm s $v^i$.  Passing to a common refinement
allows us to assume that $I=J$ and that $U_i=V_i$.  Then since $u_i$
and $v_i$ both implement $\alpha^i$ over $U_i$, an argument similar to
that above implies they differ by a central multiplier $\lambda_i:U_i\to
\hgab$.  Then it is easy to see that we get cohomologous cocycles.
One usually writes $\zeta(\alpha)$ for the class $\zeta$, and
$\zeta(\alpha)$ is called the \emph{Phillips-Raeburn obstruction}.
\begin{remark}
\label{rem-a}
When $G$ is abelian and $A$ is type~I with spectrum $X$, then
$\zeta(\alpha)$ is the classical Phillips-Raeburn obstruction of
\cite{pr2}.  That is, $\zeta(\alpha)$ is the class of the principal
$\widehat G$-bundle given by the restriction map $p:\spec{\acg} \to X$
as in \cite[Theorem~2.2]{pr2}.
\end{remark}

\begin{prop}
\label{prop-PRobs}
Let $X$ be a second countable locally compact Hausdorff space.
Suppose that $A\in\CR(X)$ and
that $\alpha:G\to\Aut(A)$ is a second countable locally compact,
locally unitary automorphism group.  Then the transition functions
\eqref{eq-*} define a class $\zeta(\alpha)$ in $H^1(X,\shgab)$ which
depends only on $(A,G,\alpha)$.  If $(A,G,\beta)$ is another such
system, then $\zeta(\alpha)=\zeta(\beta)$ if and only if $\alpha$ and
$\beta$ are exterior equivalent.  In particular, $\alpha$ is unitary
if and only if $\zeta(\alpha)=1$\footnote{We are writing the product
in $H^1$ multiplicatively; therefore $1$ denotes the trivial
element.}.
Furthermore if $A$ is stable, then every
class in $H^1(X,\shgab)$ is equal to $\zeta(\alpha)$ for some locally
unitary action $\alpha:G\to\Aut(A)$.
\end{prop}

\begin{proof}
[Proof of all but the last assertion]
We have already seen that $\zeta(\alpha)$ depends only on
$(A,G,\alpha)$.  Now suppose that $(A,G,\beta)$ is another locally
unitary action with
$\zeta(\beta)=\zeta(\alpha)$.
Then we can find a cover
$\set{U_i}_{i\in I}$ and $\ui,\vi:G\to \UM(A_{U_i})$ such that $\ui$
implements $\ai$, $\vi$ implements $\beta^i$, and such that
\begin{equation}
\label{eq-**}
u^{i}_s(x)u^{j}_s(x)^*=\gij(x)(s)=v^{i}_s(x)v^{j}_s(x)^*.
\end{equation}
Let $w^i_s(x):= u^i_s(x)v^i_s(x)^*$.  Then $s\mapsto w^i_s(\cdot)$ is
a strictly continuous map of $G$ into $\UM(A_{U_i})$.  Then it is
easy to see that $\ai$ is exterior equivalent to $\beta^i$ via $w^i$.
However, if $x\in U_{ij}$, then \eqref{eq-**} implies that
\begin{equation*}
w^i_s(x)w^j_s(x)^* = \ui_s(x)\vi_s(x)^*\vj_s(x)\uj_s(x)^* = 1.
\end{equation*}
Consequently, we can define $w_s(x)=w^i_s(x)$ if $x\in U_i$.  Since
each $w^i$ defines a strictly continuous map into $\UM(A_{U_i})$ and
since $x\mapsto\|a(x)\|$ vanishes at infinity for each $a\in A$, it is
not hard to see that $w$ is a strictly continuous map from $G$ into
$\UM(A)$.  Therefore, $\alpha$ and $\beta$ are exterior equivalent.

Conversely, if $\alpha$ and $\beta$ are exterior equivalent via
$w:G\to\UM(A)$, then with $\set{U_i}_{i\in I}$ and $\ui,\vi:G\to
\UM(A_{U_i})$ as above, we must have unimodular scalars
$\lambda_i(x)(s)$ for all $x\in U_i$ and $s\in G$ such that
$\lambda_i(x)(s) = \ui_s(x)^* w_s(x)\vi_s(x)$.  As above, we may view
these as continuous functions from $U_i$ to $\hgab$.  Also, if $x\in
U_{ij}$, then
\begin{align*}
\ui_s(x)^*\uj_s(x)& =
\ui_s(x)^*w_s(x)\vi_s(x)\vi_s(x)^*\vj_s(x)\(\uj_s(x)^*w_s(x)\vj_s(x)\)^*
\\
&= \lambda_i(x)(s)\overline{\lambda_j(x)(s)}\vi_s(x)^*\vj_s(x).
\end{align*}
It follows that $\zeta(\alpha)=\zeta(\beta)$.
\end{proof}

To prove that every class in $H^1(X,\shgab)$ arises when $A$ is
stable, we want to recall some facts about balanced tensor products.
Suppose that $A$ and $B$ are \coxalg s.  Let $I$ be the ideal in
$\atmb$ generated by
\[
\set{a\cdot f\tensor b - a\tensor f\cdot b: \text{$a\in A$, $b\in B$,
and $f\in \cox$}}.
\]
The \emph{(maximal) \cox-balanced tensor product} of $A$ and $B$ is
defined to be the quotient
\[
\atxb:= (\atmb)/I.
\]

\begin{remark}
Balanced tensor products have been studied by several authors, and
quite recently by Blanchard \cite{blanchard,blanchard2}.  In
particular if $X$ is compact, then $\atxb$ coincides with Blanchard's
$A\Xtensor^M B$.  Moreover, $\atxb$ is a \coxalg{} and, writing
$a\xtensor b$ for the image of $a\tensor b$ in $\atxb$, we have
\begin{equation}
\label{eq-(a)}
f\cdot(a\xtensor b)=f\cdot a\xtensor b = a\xtensor f\cdot
b\quad\text{for all $f\in\cox$.}
\end{equation}
We intend to discuss these and other properties of $\xtensor$
elsewhere \cite[\S2]{ew3}.
Here we will be satisfied with the special cases
outlined below.
\end{remark}

In this work, we shall always assume that $A$ and $B$ are separable,
and that $B$ is nuclear --- in fact, it will suffice to consider only
the case where $B=C_0(X,\K)$.  Then \cite[Lemma~1.1]{rw} applies and
we can identify $\Prim(\atxb)$ with
\begin{equation}
\label{eq-(b)}
\set{(P,Q)\in\Prim(A)\times\Prim(B):\sigma_A(P)=\sigma_B(Q)}.
\end{equation}
In this case, \eqref{eq-(a)} is a straightforward
consequence of \eqref{eq-(b)} and the definition of $I$.  Moreover,
for all $x\in X$,
\begin{equation}
\label{eq-fibres}
\atxb(x)\cong A(x)\tensor B(x),
\end{equation}
and $(a\xtensor b)(x)=a(x)\tensor b(x)$.
(Note that we write simply $\tensor$ when one of the factors is
nuclear.)
If $B=C_0(X,\K)$, then
$\Prim\(A\xtensor \coxk\)$ can be identified with $\Prim(A)$.
Moreover since $\coxk\cong C_0(X)\tensor \K$, the map
$a\xtensor(f\tensor T)\mapsto a\cdot f\tensor T$ extends to a
\cox-linear isomorphism of $A\xtensor\coxk$ onto $A\tensor \K$.  Notice
that if $U$ is open in $X$, then this isomorphism identifies
$(A\xtensor\coxk)_U$ with $A_U\tensor \K$.  Furthermore, if $A$ is
stable, then we can choose an isomorphism of $A\tensor \K$ and $A$
which induces the identity map on the primitive ideal spaces (assuming
$\Prim(A\tensor\K)$ has been identified with \prima)
\cite[Lemma~4.3]{pr2}.   Then $A\xtensor\coxk$ is isomorphic to $A$
and $(A\xtensor\coxk)_U$ is identified with $A_U$.

\begin{lem}[cf., {\cite[Proposition~3.10]{pr2}}]
\label{lem-PR3.10}
Suppose that $X$ is a second countable locally compact Hausdorff space
and that $A\in\CR(X)$.  Then $A\xtensor\coxk\in\CR(X)$.  If
$\alpha:G\to\Aut(A)$ and $\beta:G\to\Aut\(\coxk\)$ are
$\cox$-automorphism groups, then the diagonal action
$\alpha\tensor\beta$ on $A\tensor\coxk$ induces an action
$\alpha\xtensor\beta$ on $A\xtensor\coxk$.  If $\gamma$ is exterior
equivalent to $\alpha$ and $\delta$ is exterior equivalent to $\beta$,
then $\alpha\xtensor\beta$ is exterior equivalent to
$\gamma\xtensor\delta$.  Finally if $\alpha$ and $\beta$ are locally
unitary, then so is $\alpha\xtensor \beta$; moreover,
\[
\zeta\(\alpha\xtensor \beta)=\zeta(\alpha)\zeta(\beta)\quad
\text{in $H^1(X,\shgab)$}.
\]
\end{lem}

\begin{remark}
If $B$ is an arbitrary element of $\CR(X)$, then it seems to be
difficult to decide whether $A\xtensor B$ is in $\CR(X)$ --- even if
$A$ and $B$ are both nuclear.  However if one of the algebras is nuclear
and $A$ or $B$ has Hausdorff primitive ideal space, then one can replace
$\coxk$ by
$B$ in the above and obtain the same results.
\end{remark}

\begin{proof}
Since $\Prim(A\xtensor\coxk)$ can be identified with \prima,
$A\xtensor\coxk$ is certainly in $\CR(X)$.  If $\alpha_s$ and
$\beta_s$ are \cox-linear, then $\alpha_s\tensor\beta_s$ maps the
balancing ideal into itself
and $\alpha\xtensor\beta$ is a well-defined action on $\atxb$.

Now suppose that $u:G\to\UM(A)$ and $v:G\to\UM\(\coxk\)$ are strictly
continuous $1$-cocycles such that
$
\alpha_s(a)=u_s\gamma_s(a)u_s^*$ and $\beta_s(b)=v_s\delta_s(b)
v_s^*$ for all $s\in G$, $a\in A$, and $b\in\coxk$.
Since the image of \cox{} sits in the center of the respective
multiplier algebras, it is clear that each $u_s$ and $v_s$ commutes
with the \cox-actions.  Therefore $u_s\tensor v_s$ defines a well
defined element $w_s:=u_s\xtensor v_s$ in $\UM\(A\xtensor\coxk\)$.  The
continuity of $s\mapsto w_s(t)$ is clear for $t$ in the algebraic
tensor product $A\odot\coxk$.  This suffices to show strict continuity
as each $w_s$ has norm one.  Routine calculations show that $w_s$ is a
$1$-cocycle implementing an exterior equivalence between
$\alpha\xtensor \beta$ and $\gamma\xtensor\delta$.

Finally, suppose that $\alpha$ and $\beta$ are locally unitary.  Then
we can find a cover $\cU=\set{U_i}$ such that $\alpha^{U_i}$ is
implemented by a \hm{} $u^i:G\to\UM(A_{U_i})$ and $\beta^{U_i}$ by a
\hm{} $\vi:G\to\UM\(C_0(U_i,\K)\)$.  Let $\gamma=\set{\gij}$ and
$\eta=\set{ \eta_{ij}}$ be the corresponding cocycles representing
$\zeta(\alpha)$ and $\zeta(\beta)$.  As above, we obtain a \hm{} $w^i =
u^i\tensor_{U_i}\vi$ which implements $(\alpha\xtensor\beta)^{U_i}=
\alpha^{U_i} \tensor_{U_i} \beta^{U_i}$ on $A_{U_i}\tensor_{U_i}
C_0(U_i,\K) \cong \(A\xtensor C_0(X,\K)\)_{U_i}$.  Thus
$\alpha\xtensor\beta$ is locally unitary.  Moreover since
$w^{i}(x)= u^{i}(x)\tensor v^{i}(x)$,
\begin{align*}
w^i_s(x)&=\gij(x)(s)\uj_s(x)\tensor \eta_{ij}(x)(s)v^j_s(x) \\
&= \(\gij(x)(s)\eta_{ij}(x)(s)\)w^j_s(x).
\end{align*}
The result follows.
\end{proof}

\begin{proof}
[Proof of the final assertion in \propref{prop-PRobs}]
Let $\zeta_0\in H^1(X,\shgab)$.
As we remarked above, when $A$ is stable there is an isomorphism of
$A$ and $A\xtensor\coxk$ carrying $A_U$ onto $\(A\xtensor\coxk\)_U$.
Thus it will suffice to produce a locally unitary action $\alpha$ on
$A\xtensor\coxk$ with $\zeta(\alpha)=\zeta_0$.
It follows from  \cite[Theorem~3.8]{pr2} and Remark~\ref{rem-a}
that there is a locally
unitary action $\tilde\beta:\gab\to\Aut\(\coxk\)$ with
$\zeta(\tilde\beta)=\zeta_0$.
Now we simply lift $\tilde\beta$ to $G$.  That is $\beta_s:=
\beta_{sH}$ where $H:=\overline{[G,G]}$.
It is straightforward to check that $\zeta(\beta)=\zeta(\tilde\beta)
=\zeta_0$.
Now the result follows by applying \lemref{lem-PR3.10} to
$\alpha:=1\xtensor \beta$.
\end{proof}

\begin{remark}
Since two actions with the same Phillips-Raeburn obstruction are
exterior equivalent, the above argument makes it clear that \emph{any}
locally unitary action of $G$ on a stable \CR-algebra $A$ is lifted from
an action of $\gab$.  In fact, there is a one-to-one correspondence
between exterior equivalence classes of actions of $G$ on $A$ and
exterior equivalence classes of actions of $\gab$ on $A$.
\end{remark}

We end this section with a short discussion on locally unitary action
on continuous-trace \cs-algebras. Assume
that $A$ is a separable continuous trace $\cs$-algebra with
spectrum $X$. An action $\alpha:G\to \Aut(A)$ is called
{\em pointwise unitary} if $\alpha$ is
$C_0(X)$-linear and the action on each fibre $A(x)$ is unitary,
or, equivalently,
if each $\rho\in \hA$
can be extended to a covariant representation $(\rho, V)$
of $(A,G,\alpha)$. In general, a pointwise unitary action need not be
locally unitary (see \secref{sec-appendix}).  Despite this, pointwise
unitary actions are locally unitary under mild additional hypotheses.
The strongest result in this direction is due to Rosenberg.

\begin{thm}[{\cite[Corollary 1.2]{ros2}}]\label{thm-ros}
Let $A$ be a separable continuous-trace $\cs$-algebra with
spectrum $X$ and let $G$ be a second countable locally
compact group such that $\gab$ is compactly generated and
$H^2(G,\TT)$ is Hausdorff. Then every pointwise unitary action
of $G$ on $A$ is locally unitary.
\end{thm}

Thus as a direct corollary of this and \propref{prop-PRobs}, we obtain:

\begin{cor}\label{cor-point}
Let $A$ and $G$ be as above and assume in addition that $A$
is stable. Then the Phillips-Raeburn  obstruction map
$\alpha\to \zeta(\alpha)$ induces a bijection
between the exterior equivalence classes of pointwise unitary
actions
of $G$ on $A$ and $H^1(X,\shgab)$.
\end{cor}

%% file: repr-gr.tex
%
%

\section{Representation Groups}\label{sec-3}

In this section we want to discuss the notion of
representation groups of second countable locally compact
groups as introduced by Moore in \cite{moore2}.
Recall that if
$1\to Z\to H\to G\to 1$
is a second countable locally compact
central extension of $G$ by the abelian group
$Z$, then the transgression map
$\tg:\widehat{Z}= H^1(Z,\TT)\to H^2(G,\TT)$ is defined as follows:
Let $c:G\to H$ be a Borel section
for the quotient map $H\to G$ such that $c(eZ)=e$. Then
$\sigma(s,t):=c(s)c(t)c(st)^{-1}$
is a cocycle in $Z^2(G,Z)$ (Moore-cohomology with values
in the trivial $G$-module $Z$).
If $\chi\in \widehat{Z}$, then
$\sigma_{\chi}(s,t):=\chi(\sigma(s,t))$
defines a cocycle
$\sigma_{\chi}\in Z^2(G,\TT)$ and then
$\tg(\chi):=[\sigma_{\chi}]$
is the cohomology class of $\sigma_{\chi}$ in $H^2(G,\TT)$.

\begin{definition}[Moore]
\label{def-3.1}
Let $G$ be a second countable locally compact group and let
$H$ be a central extension of $G$ by some abelian group $Z$
such that the transgression map
$\tg:\widehat{Z}\to H^2(G,\TT)$ is bijective.
Then $H$ (or rather the extension $1\to Z\to H\to G\to 1$)
is called a {\it
representation group} for
$G$. A group $G$ is called {\em smooth} if it has a representation group
$H$.
\end{definition}

\begin{remark}\label{rem-rep}
The question of
which groups have a representation group was studied extensively
by Moore in \cite{moore2}.  If $H$ is a representation group for $G$,
then the transgression map
$\tg: \widehat{Z}\to
H^2(G,\TT)$ is actually a homeomorphism by
\cite[Theorem 2.2]{moore2} and
\cite[Theorem 6]{moore4}, so that in this case $H^2(G,\TT)$ is always
locally compact and Hausdorff. Conversely, if $G$ is almost connected
(i.e., $G/G_0$ is compact) and $H^2(G,\TT)$ is Hausdorff, then $G$ is smooth
by \cite[Proposition 2.7]{moore2}.  Since
\cite[Theorem A and following remark~(2)]{moore1} implies that
$H^2(G,\TT)$ is isomorphic to $\RR^k$ for some $k\geq 0$ whenever $G$
is a connected and simply connected Lie group, it follows that such groups
are
smooth.
If $G$ is a connected semisimple Lie group, then the universal covering group
$H$ of $G$ is a representation group for $G$ by \cite[Proposition 3.4]{moore2}.
Finally, every compact group is smooth
(see the discussion preceding
\cite[Proposition 3.1]{moore2}), and every
discrete group is smooth by \cite[Theorem~3.1]{moore2} (see also
\cite[Corollary 1.3]{para2}).
\end{remark}

We will see below that, in addition to the above, every second countable
compactly generated abelian group is smooth (\corref{corcompgen}).
To prove our results, we need the
following characterization of smooth groups.

\begin{lem}
\label{cocycle}
Let $G$ be a second countable locally compact group.
Then $G$ is smooth if and only if there exists a
second countable locally compact Hausdorff topology on $H^2(G,\TT)$
and a Borel cocycle $\zeta\in Z^2\(G, \specnp{H^2(G,\TT)}\)$ such
that for each $[\om]\in H^2(G,\TT)$ evaluation of $\zeta$ at
$[\om]$ gives an element in $Z^2(G,\TT)$ representing $[\om]$.
\end{lem}
\begin{proof}
Assume $1\to Z\to H\to G\to 1$ is a
representation group for
$G$.
Since $\tg$ is an isomorphism, we can define an isomorphism
$\tg_*:Z\to \specnp{H^2(G,\TT)}$ by $ \tg_*(z)([\om])=\tg^{-1}([\om])(z)$.
Let $c:G\to H$ be a Borel cross section with
$c(e_G)=e_H$ ($e_G$ denoting the unit in $G$).
For $s,t\in G$, define $\zeta(s,t):=\tg_*(c(s)c(t)c(st)^{-1})$. Then
$\zeta\in Z^2(G, \specnp{H^2(G,\TT)})$, and if
$[\om]\in H^2(G,\TT)$, then
we obtain $\zeta(s,t)([\om])=\tg^{-1}([\om])(c(s)c(t)c(st)^{-1})$.
Thus, by the definition of $\tg$ (see the discussion above)
$(s,t)\mapsto \zeta(s,t)([\om])$ is a cocycle representing $[\om]$.
This is what we wanted.

For the converse, assume that there is a second countable
locally compact Hausdorff topology on $H^2(G,\TT)$ and
let $\zeta\in Z^2(G,\specnp{H^2(G,\TT)})$ be such that evaluation at
each $[\om]\in H^2(G,\TT)$ gives a representative for $[\om]$.
Let $H$ denote the extension $G\times_\zeta \specnp{H^2(G,\TT)}$ of $G$
by $\specnp{H^2(G,\TT)}$ given by
$\zeta$. This is the set $G\times \specnp{H^2(G,\TT)}$
with multiplication defined by
$$(s,\chi)(t,\mu)=(st, \zeta(s,t)\chi\mu),$$
and equipped with the unique locally compact group topology inducing
 the product Borel structure on $G\times \specnp{H^2(G,\TT)}$
(see \cite[Theorem 7.1]{m5}).
Then $$1\to \specnp{H^2(G,\T)}\to H\to G\to 1$$ is a central
extension of $G$ by $\specnp{H^2(G,\T)}$.
Let $c:G\to H$ denote the canonical section $c(s)= (s, 1)$.
Then the transgression map is a map from $H^2(G,\TT)\cong
\specnp{(\specnp{H^2(G,\T)})}$
to itself, and for each $[\om]\in H^2(G,\T)$ we obtain a representative
for $\tg([\om])$ by taking the cocycle
\[\nu(s,t)=(c(s)c(t)c(st)^{-1})([\om])=\zeta(s,t)([\om]).\]
Therefore the transgression map $\tg:H^2(G,\TT)\to H^2(G,\TT)$
is the identity, and $H$ is a representation group for $G$ as required.
\end{proof}

\begin{remark}\label{rem-isom}
If $1\to Z\to H\to G\to 1$ is a representation group
for $G$ and if $\zeta\in Z^2(G,\specnp{H^2(G,\TT)})$ is
constructed as above such that the evaluation map is
the identity on $H^2(G,\TT)$, then $G\times_{\zeta}\specnp{H^2(G,\TT)}$ is
actually isomorphic $H$.
An isomorphism is given by $(s,[\sigma])\mapsto c(s)(\tg_*)^{-1}([\sigma])$,
where $c:G\to H$ denotes the Borel section defining $\zeta$ as above.
\end{remark}

We now show that, under some weak additional
assumptions, the direct product of
smooth groups is again smooth.

\begin{prop}\label{represent}
Suppose that $G_1$ and $G_2$ are smooth and let
$B(G_1,G_2)$ denote the group of continuous bicharacters $\chi:G_1\times G_2\to
\TT$. If $B(G_1,G_2)$ is locally compact with respect to the compact open
topology, then $G_1\times G_2$ is smooth. In particular, if
the abelianizations $(G_1)_{\subab}$, $(G_2)_{\subab}$ of $G_1$
and $G_2$ are compactly generated, then $G_1\times G_2$ is smooth.
\end{prop}
\begin{proof}
In fact we are going to construct a representation group for
$G_1\times G_2$ as follows: Choose central extensions
$$1\to \specnp{H^2(G_i,\TT)}\to H_i\stackrel{q_i}{\to}G_i\to 1$$
for $i=1,2$ such that the respective transgression maps are both equal
to the identity
(see Lemma~\ref{cocycle} and Remark~\ref{rem-isom}).
By assumption $B(G_1,G_2)$ is locally compact, so the dual group
$\specnp{B(G_1,G_2)}$ is also a locally compact group.
For each pair $(s_1,s_2)\in G_1\times G_2$ define $\eta(s_1,s_2)\in
\specnp{B(G_1,G_2)}$ by $\eta(s_1,s_2)(\chi)=\chi(s_1,s_2)$,
$\chi\in B(G_1,G_2)$.
Let
$$H=H_1\times H_2\times \specnp{B(G_1,G_2)}$$
with multiplication defined by
\begin{equation}
\label{eq-extra}
(h_1, h_2, \mu)(l_1, l_2,\nu)=
\big(h_1l_1, h_2l_2,\mu\nu\eta(q_1(h_1),q_2(l_2))\big).
\end{equation}
Then clearly
$$Z= \specnp{H^2(G_1,\TT)}\times
\specnp{H^2(G_2,\TT)}\times\specnp{B(G_1,G_2)}$$
is a central subgroup of  $H$, and we obtain a short exact sequence
$$1\to Z\to H\to G_1\times G_2\to 1.$$
We claim that $H$ is a representation group for $G_1\times G_2$.

For this recall that if
$\om_1\in Z^2(G_1,\TT)$, $\om_2\in Z^2(G_2,\TT)$ and $\chi\in B(G_1,G_2)$,
then $\om_1\otimes\om_2\otimes\chi$
defined by
\[
\om_1\otimes\om_2\otimes\chi((s_1,s_2),(t_1,t_2))=
\om_1(s_1,t_1)\om_2(s_2,t_2)\chi(s_1,t_2).
\]
is a cocycle in $Z^2(G_1\times G_2,\TT)$.
By \cite[Theorem~9.6]{m4} and
\cite[Propositions~1.4 and 1.6]{klepp} we
know that the map
$$([\om_1],[\om_2],\chi)\mapsto [\om_1\otimes\om_2\otimes\chi]$$
is an (algebraic) isomorphism of $H^2(G_1,\TT)\times H^2(G_2,\TT)\times
B(G_1,G_2)$ onto $H^2(G_1\times G_2,\TT)$, from which it follows that
$Z$ is isomorphic to $\specnp{H^2(G_1\times G_2,\TT)}$.
Now choose Borel sections $c_i:G_i\to H_i$ and define
$\zeta_i\in Z^2(G_i, \specnp{H^2(G_i,\TT)})$ by
$\zeta_i(s,t)=c_i(s)c_i(t)c_i(st)^{-1}$, for $s,t\in G_i$.
Since the transgression maps are both equal to the identity map, we see that
evaluation of $\zeta_i$ at $[\om_i]\in H^2(G_i,\TT)$ is a cocycle
representing $[\om_i]$. Defining  $c:G_1\times G_2\to H_1\times H_2\times
\specnp{B(G_1,G_2)}$ by $c(s_1,s_2)=(c_1(s_1), c_2(s_2), 1)$, we easily compute
$$
\zeta((s_1,s_2),(t_1,t_2)):=c(s_1,s_2)c(t_1,t_2)c(s_1t_1, s_2t_2)^{-1}=
(\zeta_1(s_1,t_1),\zeta_2(s_2,t_2),
\eta(s_1,t_2)).$$
Thus, evaluating $\zeta$ at
$[\om_1\otimes\om_2\otimes\chi]\in H^2(G_1\times G_2,\TT)$ gives a
cocycle representing this class. This proves the claim.
The final assertion now follows from the fact that
$B(G_1, G_2)=B((G_1)_{\subab},(G_2)_{\subab})$ and
\cite[Theorem 2.1]{klepp}.
\end{proof}

\begin{cor}\label{corcompgen}
Every second countable compactly generated abelian group is smooth.
\end{cor}
\begin{proof}
By the structure theorem for compactly generated abelian groups
(\cite[Theorem~9.8]{hr}), we know
that $G\cong \RR^n\times K\times \ZZ^m$ for some $n,m\geq 0$
and some compact group $K$. By the results mentioned in Remark~\ref{rem-rep},
it follows that $\RR^n,K$ and $\ZZ^m$ are smooth.
Now apply the proposition.
\end{proof}

The example given at the bottom of \cite[p.\ 85]{moore2} shows that
there are nonsmooth abelian locally compact groups.
The group constructed there is a direct product of $\RR$
with an infinite direct sum of copies of $\ZZ$. Thus it also provides
an example of two smooth groups  whose direct product is not
smooth.  Thus the assumption on $B(G_1,G_2)$ in
Proposition~\ref{represent}
is certainly not superfluous.

It is certainly interesting to see specific examples of
representation groups. Some explicit constructions can be found in
\cite{moore2} and
\cite[Corollary 1.3 and Examples 1.4]{para2}. For instance,
if $G=\RR^2$ (resp.\ $G=\ZZ^2$), then the three
dimensional Heisenberg group (resp.\ discrete Heisenberg group)
is a representation group for $G$. In the following example,
we use Proposition~\ref{represent} to construct representation groups
for $\RR^n$.

\begin{example}\label{ex-rn}
Let $G=\RR^n$ and, as a set, let $H_n=\RR^{{n(n+1)}/{2}}$.
We write an element of $H_n$ as
$\mathfrak s= (s_i, s_{j,k})$, $1\leq i\leq n$, ${1\leq j<k\leq n}$.
Define
multiplication on $H_n$ by
$\mathfrak s\mathfrak t=((st)_i, (st)_{j,k})$ with
$$(st)_i:=s_i+t_i\quad\text{and}\quad
(st)_{j,k}:=s_{j,k}+t_{j,k}+s_jt_k.$$
Then $H_n$ is clearly a central extension of $\RR^n$ by
$\RR^{{(n-1)n}/{2}}$.

We claim that $H_n$ is a representation group for $\RR^n$.
Since $H^2(\RR,\TT)$ is trivial, this is certainly true for $n=1$.
For the step $n\to n+1$ assume that $H_n$ is a representation group for
$\RR^n$. For $\vec s=(s_1,\ldots,s_n)\in \RR^n$ define $\chi_{\vs}
\in B(\RR^n,\RR)$ by
$$\chi_\vs((t_1,\ldots, t_n), r):=e^{ir(s_1t_1+\cdots+s_nt_n)}.$$
Since $\chi_\vs(\vec t,r)=\chi_{r\vs}(\vec t,1)$,
$s\mapsto \chi_{\vec s}$ 
is an isomorphism of $\RR^n$ onto $B(\RR^n,\RR)$, and we
see that $\specnp{B(\RR^n,\RR)}$ is isomorphic to $\RR^n$ via
the map $\vec t\mapsto \eta_{\vec t}$
defined by $\eta_{\vec t}(\chi_\vs)=\chi_\vs(\vec t,1)$, for $\vec t\in \RR^n$.
Moreover, if $(\vs,r)\in \RR^n\times \RR$, and
$\eta(\vs,r)\in \specnp{B(\RR^n,\RR)}$ is defined by
$\eta(\vs,r)(\chi_{\vec t})=\chi_{\vec t}(\vs,r)$,
then we get the identity $\eta(\vs,r)=\eta_{r\vs}$.
It follows now from  Proposition~\ref{represent} and \eqref{eq-extra}
that
$H'=H_n\times \RR\times \RR^n$ with multiplication defined by
$$\big((s_i, s_{j,k}), r, \vec t\big)\big((s_i', s_{j,k}'), r', \vec t'\big)=
\big((s_i, s_{j,k})(s_i', s_{j,k}'), r+r',\vec t+\vec t'+r'\vs\big)$$
is a representation group for $\RR^{n+1}$. Putting $s_{n+1}=r$
and $s_{j,n+1}=t_j$, $1\leq j\leq n$, we see that this formula coincides with
the multiplication formula for $H_{n+1}$.
\end{example}

Using similar arguments, it is not hard to show that a representation group for
$\ZZ^n$ is given by the integer subgroup of $H_n$ constructed above,
i.e., assuming that all $s_i$ and $s_{j,k}$ are integers.
Notice that the group $H_n$ constructed above is isomorphic to the
connected and simply connected two-step nilpotent Lie group
corresponding to the universal two-step nilpotent
Lie algebra generated  by $X_1, \ldots, X_n$ and the commutators $[X_j,X_k]$,
$1\leq j<k\leq n$. Note that any connected two-step
nilpotent Lie group is a quotient of one of these
groups (e.g., see \cite[p.~409]{brown}).

We conclude this section with a discussion of which conditions will
imply that all representation groups of a given group $G$ are
isomorphic. Schur observed that even finite groups can have
nonisomorphic representation groups \cite{schur}, and Moore considers
the case for $G$ compact or discrete in\cite[\S 3]{moore2}.  Here we
give a sufficient condition for the uniqueness of the representation
group (up to isomorphism) valid for all smooth $G$.

\begin{prop}\label{prop-unique}
Let $G$ be smooth and let $Z:=\specnp{H^2(G,\TT)}$.
Then the representation groups
of $G$ are unique (up to isomorphism of groups) if
every abelian extension $1\to Z\to H\to \gab\to 1$
splits.
In particular,
if
$\gab$ is isomorphic to
$\RR^n\times
\ZZ^m$  or if $Z$ is isomorphic to $\RR^n\times \TT^m$, for some $n,m\geq 0$,
then all representation groups of $G$ are isomorphic.
\end{prop}
\begin{proof}
Let $1\to Z\to H_1\to G\to 1$ and $1\to Z\to H_2\to G\to 1$
be two representation groups of $G$.
By Lemma~\ref{cocycle} and Remark~\ref{rem-isom} we
may assume that (up to isomorphism) both extensions are given by
cocycles $\zeta_1,\zeta_2\in Z^2(G,Z)$ such that the
transgression maps $H^2(G,\TT)=\widehat{Z}\to H^2(G,\TT)$ induced by
$\zeta_1$ and $\zeta_2$ are the identity maps.
Now let $\sigma=\zeta_1\circ \zeta_2^{-1}\in Z^2(G,Z)$ and let
$1\to  Z\to L\to G\to 1$ denote the extension defined by $\sigma$.
We want to show that this extension splits (then $\sigma\in B^2(G,Z)$
and $[\zeta_1]=[\zeta_2]\in H^2(G,Z)$).

Since $\chi\circ \sigma=(\chi\circ \zeta_1)\cdot (\chi\circ \zeta_2^{-1})$
and $[\chi\circ \zeta_1]=[\chi\circ \zeta_2]\in H^2(G,\TT)$, it follows
that the transgression map $\widehat{Z}\to H^2(G,\TT)$ induced by
$\sigma$ is trivial. But this implies that
any character of $Z$ can be extended to a character of $L$,
which implies that
$\widehat{L}_{\subab}$ is an extension
$1\to \hgab\to \widehat{L}_{\subab}\to \widehat{Z}\to 1$.
By assumption (using duality), this extension splits.
Thus we find an injective homomorphism
$\chi\mapsto \mu_{\chi}$ from $\widehat{Z}\to\widehat{L}_{\subab}$
such that
each $\mu_{\chi}$ is an extension of $\chi$ to $L$.
Let $\tilde{G}=\{s\in L: \mu_{\chi}(s)=1\;\text{for all}\;\chi\in
\widehat{Z}\}$.
Then $\tilde{G}\cap Z=\{e\}$ and $\tilde{G}\cdot Z =L$.
To see the latter, let $l\in L$ and let $z\in Z$ such that
$\mu_{\chi}(l)=\chi(z)$ for all $\chi\in \widehat{Z}$. Then
$lz^{-1}\in \tilde{G}$. It follows that the quotient map $q:L\to G$
restricts to an isomorphism $\tilde{G}\to G$.  This proves all but the
final statement.

By duality, it suffices to prove the final assertion only for $\gab$
isomorphic to $\R^n\times\Z^m$.  By induction, it suffices to consider
only the cases $\Z$ and $\R$.  Since the first is straightforward,
we will show only that if $G$ is abelian, then any continuous open
surjection $q:G\to\R$ has a continuous section.  By
\cite[Theorem~24.30]{hr}, we may assume that $G=\R^m\times H$, where
$H$ has a compact, open subgroup.  It follows that $q\restr{\R^m}$ is
surjective.  Thus there is an $x\in\R^m$ such that $q(x)=1$.  Then we
can define $q^*:\R\to\R^m$ by $q(\lambda)=\lambda\cdot x$.
\end{proof}

%% file: brauer.tex
%
%

\section{The Brauer Group for Trivial Actions}\label{sec-4}

In this section we want to give a precise description
of the set $\E_G(X)$
of exterior equivalence classes of $C_0(X)$-linear
actions $\beta:G\to \Aut\(C_0(X,\K)\)$ in the case where $G$ is
smooth,
$\gab$ is compactly
generated, and $X$ is a second countable locally compact
Hausdorff space.  This analysis also allows a
description of the Brauer group $\Br_G(X)$ of
\cite{ckrw}
for a trivial $G$-space $X$,
and, more generally, a description
of the set of exterior equivalence classes of
locally inner actions $\alpha:G\to \Aut(A)$
when $A\in \CR(X)$.

A \cs-dynamical system $(A,G,\alpha)$ is called a
\emph{$C_0(X)$-system} if $A$ is a \coxalg{} and each $\alpha_s$ is
\cox-linear.  Two systems $(A,G,\alpha)$ and $(B,G,\beta)$ are
Morita equivalent if there is a pair $(\X,\mu)$ consisting of an $A
\sme B$-\ib{} $\X$ and a strongly continuous action $\mu$ of $G$ on
$\X$ by linear transformations such that
\[
\alpha_s\(\lip A<x,y>\)=
\blip A<\mu_s(x),\mu_s(y)>
\quad\text{and}\quad
\beta_s\(\rip B<x,y>\)=
\brip B<\mu_s(x),\mu_s(y)>
\]
 for all $x,y\in
\X$ and $ s\in G$.
The actions of $A$ and $B$ on $\X$ extend to the multiplier algebras
$M(A)$ and $M(B)$.  In particular, if $A$ and $B$ are \coxalg s, then
$\X$ is both a left and a right \cox-module.  We say that two
\cox-systems $(A,G,\alpha)$ and $(B,G,\beta)$ are \emph{\cox-Morita
equivalent} if they are Morita equivalent and if $f\cdot x=x\cdot f$
for all $x\in\X$ and $f\in\cox$.
If we let $G$ act trivially on $X$, then the equivariant Brauer group $\brg(X)$
of \cite{ckrw} is the collection of $\cox$-Morita equivalence classes
of $\cox$-systems $(A,G,\alpha)$ where $A$
is a separable continuous-trace $\cs$-algebra with spectrum
$X$.
Then $\Br_G(X)$ forms an
abelian group \cite[Theorem~3.6]{ckrw}.
Recall that the group multiplication is defined using the
balanced tensor product
\begin{equation}
\label{eq-prod*}
[A,\alpha][B,\beta] = [A\Xtensor B,\alpha\Xtensor \beta].
\end{equation}
The identity is the class of $\(\cox,\id\)$ and the inverse of
$[A,\alpha]$ is given by the class of the conjugate system $(\overline
A,\bar\alpha)$.

The collection of $[A,\alpha]$ in $\brg(X)$ such that the
Dixmier-Douady class of $A$ is zero is a subgroup.  Note that each such
element has a representative of the form $\(\coxk,\alpha\)$.
That we can identify this subgroup with
\egx{} follows from the next proposition.  In particular, \egx{} is also an
abelian
group with multiplication given by \eqref{eq-prod*} after identifying
$\coxk\xtensor\coxk$ with \cox{} via a \cox-isomorphism.

\begin{prop}
\label{prop-help}
Suppose that $\alpha,\gamma:G\to\Aut\(\coxk\)$ are \cox-actions such
that $[\coxk,\alpha]=[\coxk,\gamma]$ in $\brg(X)$.  Then $\alpha$ and
$\gamma$ are exterior equivalent.
\end{prop}

\begin{proof}
Since $\coxk$ is stable, it follows from \cite[Lemma~3.1]{ckrw} that
$\gamma$ is exterior equivalent to an action $\beta$ of the form
$\beta=\Phi\circ \alpha$ for some \cox-linear automorphism $\Phi$ of
$\cox$.  Thus it will suffice to see that $\beta$ is exterior
equivalent to~$\alpha$.

Since a $C_0(X)$-linear automorphism
of $C_0(X,\K)$ is locally inner \cite{pr1},
we may find an open  cover $\set{U_i}_{i\in I}$
of $X$ and continuous functions $u^i$ from $U_i$ to $\U(\H)$
for each $i\in I$ such that
$\Phi(f)(x)=u_i(x)f(x)u_i(x)^*$ for all $f\in C_0(U_i,\K)$.
Moreover, on each overlap $U_{ij}$, there exist continuous functions
$\chi_{ij}\in C(U_{ij},\TT)$
such that $u_j(x)=\chi_{ij}(x)u_i(x)$ for all $x\in U_{ij}$.
Since $\beta_s$ is a \cox-automorphism, there are automorphisms
$\beta^x_s$ for each $x\in X$ such that $\beta_s(f)(x)=\beta^x_s\(f(x)\)$.
Since $\alpha_s=\Ad\Phi\circ \beta_s$ it follows for all $x\in U_i$,
\[
\alpha_s(f)(x)=u_i(x)^*\bar\beta^x_s\(u_i(x)\)
\beta^x_s(f)(x)\bar\beta^x_s(u_i(x)^*)u
_i(x),
\]
where $\bar\beta_s^x$ is the canonical extension of $\beta_s^x$ to
$M\(A(x)\)$.
If $x\in U_{ij}$, then
\[
u_j(x)^*\beta^x_s\(u_j(x)\)=\overline{\chi_{ij}(x)}u_i(x)^*\bar\beta^x_s
\({\chi_{i
j}(x)}u_i(x)\)=
u_i(x)^*\beta^x_s\(u_i(x)\).
\]
Consequently, we can define a map from $G$ to $U(\H)$
by $v_s(x)=u_i(x)^*\bar\beta^x_s\(u_i(x)\)$.
Moreover, $s\mapsto v_s$ is strictly continuous, and  we have
$\alpha=\Ad v\circ \beta$.  Thus we only need to verify that $v$ is a
$1$-cocycle.
For all $x\in
X$ we get
\begin{align*}
v_{st}(x)&=u_i(x)^*\bar\beta^x_{st}\(u_i(x)\)=u_i(x)^*\bar
\beta^x_s\(\bar\beta^x_t\(u_i(x)\)\)
\\
&=u_i(x)^*\bar\beta^x_s\(u_i(x)\)\bar\beta^x_s\(u_i(x)^*\bar
\beta^x_t\(u_i(x)\)\)=v_s(x)\bar\beta
^x_s\(v_t(x)\),
\end{align*}
which implies $v_{st}=v_s\beta_s(v_t)$.
\end{proof}

\begin{remark}
\label{rem-help}
Suppose that $\gamma:G\to \Aut\(\coxk\)$ is a \cox-automorphism group.
In the sequel, we will write $\igamma$ for the ``inverse''
automorphism group in $\egx$.  That is,
$[\coxk,\gamma]^{-1}:=[\coxk,\igamma]$.  \propref{prop-help} implies
that $\igamma$ is unique up to exterior equivalence and that
$\igamma\Xtensor\gamma$ is exterior equivalent to $\id\Xtensor\id$.
\end{remark}

The next lemma is a mild strenthening of \cite[Lemma~3.3]{doir} to our
setting.

\begin{lem}
\label{lem-fix}
Suppose that $\beta:G\to\Aut\(\coxk\)$ is a \cox-linear action and
that $[\omega_x]$ is the Mackey obstruction for the induced
automorphism group $\beta^x$ on the fibre over $x$.  Then the
\emph{Mackey obstruction map} $\phi^\beta:X\to H^2(G,\T)$ given by
$\phi^\beta(x):=[\omega_x]$ is continuous.
\end{lem}

\begin{proof}
Fix $x_0\in X$ and suppose that $\set{x_n}$ is a sequence converging
to $x_0$ in $X$.  It will suffice to show that $[\omega_{x_n}]$
converges to $[\omega_{x_0}]$ in $H^2(G,\T)$.  Let $M$ be the compact
set $\set{x_n}_{n=1}^\infty\cup \set{x_0}$.  Let $\beta^M$ be the
induced action on $C(M,\K)$.  Since $\beta$ and $\beta^M$ induce the
same action on the fibres, we have $\phi^{\beta^M}=\phi^\beta\restr
M$.  Thus it will suffice to see that the former is continuous.  But
$H^2(M;\Z)$ is trivial; any principal $\T$-bundle over $M$ is locally
trivial and therefore trivial.  If follows from the Phillips-Raeburn
exact sequence \cite[Theorem~2.1]{pr1}
that $\beta^M$ is inner.  As in Remark~\ref{rem-fix}, there is an
obstruction to $\beta^M$ being unitary given by a cocycle $\zeta\in
Z^2\(G,C(M,\T)\)$, and $\phi^{\beta^M}(x)=[\zeta(x)]$.  Since for each
$s,t\in G$, $\zeta(s,t)$ is continuous, it follows that $\zeta(x_n)$
converges to $\zeta(x)$ pointwise.  Therefore $\zeta(x_n)\to\zeta(x)$
in $Z^2(G,\T)$ \cite[Proposition~6]{moore3}.  Since $H^2$ has the
quotient topology, the result follows.
\end{proof}

Using the above lemma, the discussion in the introduction shows
that there is a \hm{} $\Phi$ from
\egx{} to $\cxhtt$ which assigns to each $[\alpha]$ in \egx{} its
``Mackey obstruction map'' $\phi^\alpha$.

\begin{thm}
\label{brauer}
Suppose that $G$ is smooth.  Then
the homomorphism
$\Phi: \E_G(X)\to C\(X,H^2(G,\TT)\)$ given by $[\beta]\mapsto
\varphi^{\beta}$ is surjective and the short exact sequence
\[
1 \arrow{e} \ker\Phi\arrow{e} \E_G(X)\arrow{e}
C\(X,H^2(G,\TT)\)\arrow{e} 1
\]
splits. If, in addition, $\gab$ is compactly generated, then
\[
\E_G(X)\cong H^1(X,\shgab)\oplus
C\(X,H^2(G,\TT)\)\]
as abelian groups.
\end{thm}
\begin{proof}
We have to construct a splitting homomorphism
$\Phi^*:C\(X, H^2(G,\TT)\)\to \E_G(X)$ for $\Phi$.
Recall from \cite[Theorem 5.1]{ckrw} and
\cite[Proposition 3.1]{horr} that there is a canonical
homomorphism
$\mu:H^2\(G,C(X,\TT)\)\to
\E_G(X)$ defined as follows: Let $\sigma\in Z^2\(G,
C(X,\TT)\)$, and let $L^{\sigma(x)}$ denote
the left regular $\sigma(x)$-representation, where
$\sigma(x)$ denotes evaluation of $\sigma$ at $x\in X$.
A representative for the class $\mu([\sigma])\in \mathcal
E_G(X)$ is then given by the action $\beta^{\sigma}:G\to
\Aut\(C_0(X,\K\(L^2(G)\)\)$ defined by
$$\beta^\sigma_s(f)(x)=\Ad L^{\sigma(x)}_s\(f(x)\),\quad f\in
C_0(X,\K(L^2(G))).$$
(If $G$ is finite
we have to stabilize this action in order to get an action on
$C_0\(X,\K\(L^2(G)\)\)\otimes \K\cong C_0(X,\K)$.)

By Lemma~\ref{cocycle} we know that
$H^2(G,\TT)$ is  locally compact and that
there exists an element $\zeta\in
Z^2(G,\specnp{H^2(G,\TT)})$  such that evaluation of $\zeta$
at a point
$[\om]\in H^2(G,\TT)$ is a cocycle representing $\om$.
If we give $C(H^2(G,\TT),\TT)$ the compact-open topology, then we can
view $\specnp{H^2(G,\TT)}$ as a subset of $C(H^2(G,\TT),\TT)$.
Furthermore, if $\varphi\in C\(X,H^2(G,\TT)\)$, then
$\zeta\circ \varphi(s,t)(x):=\zeta(s,t)\(\varphi(x)\)$ defines a
Borel cocycle $\zeta\circ \varphi\in Z^2\(G,C(X,\TT)\)$.

We claim that
$\Phi^*(\varphi):=\mu([\zeta\circ \varphi])$
defines a splitting
homomorphism  for $\Phi$. To see that it is a homomorphism
just notice that if $\varphi,\psi\in C\(X,H^2(G,\TT)\)$, then since
$\zeta(s,t)\in \specnp{H^2(G,\TT)}$,
$\zeta(s,t)\(\varphi(x)\)\zeta(s,t)\(\psi(x)\)=
\zeta(s,t)\(\varphi(x)\psi(x)\)$.
 Thus $\varphi\mapsto
[\zeta\circ \varphi]$ is a homomorphism of $C\(X,
H^2(G,\TT)\)$ into $H^2\(G, C(X,\TT)\)$.
By the construction of
$\mu$ we can choose a representative $\beta$ for
$\mu([\zeta\circ\varphi])$ such that $\beta^x$ is implemented
by a
$\zeta\(\varphi(x)\)$-representation $V:G\to \U(\H)$.
Since $\zeta\(\varphi(x)\)$ is a representative for
$\varphi(x)$, it follows that $\Phi\circ \Phi^*=\id$.

We have shown that if $G$ is smooth, then
$$1\arrow{e} \ker\Phi\arrow{e} \E_G(X)\arrow{e}
C\(X,H^2(G,\TT)\)\arrow{e} 1$$
is a split short exact sequence. If, in addition,
$\gab$ is compactly generated, then we know from
Corollary~\ref{cor-point} that the
Phillips-Raeburn obstruction
$\beta\mapsto \zeta(\beta)\in H^1(X,\shgab)$ of Proposition
\ref{prop-PRobs} defines a bijection of $\ker\Phi$ onto
$H^1(X,\shgab)$, which by Lemma~\ref{lem-PR3.10} is
multiplicative. This completes the proof.
\end{proof}

Let $F:\Br_G(X)\to H^3(X,\ZZ)$ denote the
forgetful homomorphism described in the introduction.
It admits a natural splitting map,
which assigns to an element
$\delta\in H^3(X,\ZZ)$
the (equivalence class of the) system $(A_\delta,G,\id)$,
where $A_\delta$ is
the unique stable continuous-trace $\cs$-algebra with
Dixmier-Douady
invariant $\delta$ and $\id$ denotes the trivial action of
$G$ on $A_\delta$. Since $\ker F$ is naturally isomorphic
to $\mathcal E_G(X)$ by Proposition \ref{prop-help}, we
obtain the following as an immediate corollary.

\begin{cor}\label{cor-brauer}
Suppose that $G$ is smooth and that $\gab$ is
compactly generated. Then, for any trivial $G$-space $X$, we have a
group isomorphism
$$\brg(X)\cong
H^1(X,\shgab)\oplus C\(X,H^2(G,\TT)\)\oplus H^3(X;\ZZ),$$
where $H^3(X;\ZZ)$ denotes third integral \v Cech cohomolgy.
\end{cor}

We conclude this section with a discussion of some
special cases.

\begin{example}
\label{ex-brauer}
 If $G$ is connected and $H^2(X,\ZZ)$ is countable
then it follows from \cite[\S6.3]{ckrw} that the
homomorphism $\mu:H^2\(G,C(X,\TT)\)\to \E_G(X)$
described in the proof of Theorem~\ref{brauer} is
actually an isomorphism (in particular, all $C_0(X)$-actions are
inner). Under this isomorphism, the  Mackey obstruction map
$\Phi:
\E_G(X)\to C\(X, H^2(G,\TT)\)$ corresponds to the
evaluation map
$H^2(G,C(X,\TT)\to C\(X, H^2(G,\TT)\)$ and the kernel
of $\Phi$ corresponds to the subgroup $H^2_{\text{pt}}\(G, C(X,\TT)\)$
of pointwise trivial elements in $H^2\(G, C(X,\TT)\)$.
\end{example}

If $G$ is not connected, there are usually lots of
$C_0(X)$-linear  actions of $G$ on
$C_0(X,\K)$ which are not inner, for instance, if $G=\ZZ^n$,
then $H^2\(\ZZ^n, C(X,\TT)\)\cong C\(X, H^2(\ZZ^n,\TT)\)$ by
\cite[Corollary 1.5]{judymc},  but $H^1(X,\widehat{\ZZ}^n)$
is often nontrivial (for instance for $G=\ZZ$ and
$X=S^2$). In any case, if $G$ is smooth and if $\zeta$ is as in
\lemref{cocycle},  then the map
$\varphi\mapsto
\zeta\circ \varphi$ from $ C\(X, H^2(G,\TT)\)$ to $H^2\(G, C(X,\TT)\)$ is a
splitting homorphism for the exact sequence
$$1\arrow{e} H_{\text{pt}}^2\(G, C(X,\TT)\)\arrow{e} H^2\(G,
C(X,\TT)\)\arrow{e}
C\(X, H^2(G,\TT)\)\arrow{e} 1.$$

\begin{example}
\label{ex-two}
If $G$ is smooth and $\gab$ is a vector group
(i.e., $\gab$ is isomorphic to some $\RR^l$ for $l\geq 0$)
then $H^1(X,\shgab)=0$ and $\E_G(X)\cong
C\(X,H^2(G,\TT)\)$. This applies to all simply connected
and connected Lie groups.
\end{example}

\begin{example}
If $G$ is any second countable
locally compact group such that $H^2(G,\TT)$ is trivial
(e.g., if $G=\RR$, $\TT$, $\ZZ$, or any connected and simply
connected semisimple Lie group), then
$G$ serves as a representation group for itself.
Hence, if
$\gab$ is also compactly generated, then $\E_G(X)\cong
H^1(X,\shgab)$. If, in addition,
$\gab$ is  a vector group, then it follows from Example~\ref{ex-two} that
$\E_G(X)$ is trivial.
\end{example}
\begin{example}
It follows from the previous example, that if
$G$ is any connected and simply connected
Lie group with $H^2(G,\TT)$ trivial, then $\mathcal
E_G(X)$ is trivial. Since $H^2\(G, C(X,\TT)\)$ imbeds injectively
into $\E_G(X)$ by \cite[\S6.3]{ckrw},
it follows that for such groups
$H^2\(G, C(X,\TT)\)$ is trivial for all $X$. For compact $X$
this was shown in  \cite[Theorem 2.6]{heros}.
\end{example}

%% file: brauer1.tex
%
%

\section{Locally Inner Actions on \CR-Algebras}\label{sec-5}

 In this section we want to use our
description  of $\E_G(X)$ to describe
locally inner actions on elements of $\CR(X)$.
The next lemma provides an
analogue for the Mackey-obstruction map $\beta\mapsto
\varphi^{\beta}\in C\(X,H^2(G,\TT)\)$ of Theorem~\ref{brauer} in
case of locally inner actions on general elements of $\CR(X)$.

\begin{lem}\label{lem-function}
Let $G$ be a second countable locally compact group,
$A\in \CR(X)$, and let $\alpha:G\to \Aut(A)$ be locally
inner. For each $x\in X$ let $U$ be an open  neighborhood
of $x$ such that the restriction
$\alpha^{U}:G\to\Aut(A_{U})$ of $\alpha$ is inner, and let
$[\sigma]\in H^2\(G, C(U,\TT)\)$ be the obstruction for
$\alpha^{U}$ being unitary. Then
$\varphi^{\alpha}(x):=[\sigma(x)]$
determines a well defined continuous map
$\varphi^{\alpha}:X\to H^2(G,\TT)$.  If $\beta$ is exterior equivalent
to $\alpha$, then $\varphi^\alpha=\varphi^\beta$.
\end{lem}
\begin{proof}
We have to show that if $U_1$ and $U_2$ are two
open neighborhoods of $x$ such that, for $i=1,2$,
there exist cocycles $\sigma^i\in Z^2\(G, C(U_i,\TT)\)$ and
$\sigma^i$-homomorphism $V^i:G\to \UM\(C_0(U_i,\K)\)$ which
implement $\alpha^i:=\alpha^{U_i}$, then
$[\sigma^1(x)]=[\sigma^2(x)]$ in $H^2(G,\TT)$.
Let $U_{ij}:=U_1\cap U_2$ and let $V^{ij}_s$ denote the image
of $V^i_s$ in $\UM(A_{U_{ij}})$, $i,j\in \{1,2\}$, and let
$\sigma^{ij}$ denote the restriction of $\sigma^i$ to
$U_{ij}$; that is,
$\sigma^{ij}(s,t)(x)=\sigma^i(s,t)(x)$ for all $x\in U_1\cap
U_2$. Then $V^{ij}$ is a $\sigma^{ij}$-homomorphism which
implements the restriction
$ \alpha^{U_{ij}}:G\to \Aut(A_{U_{ij}})$.
Thus it follows that
$[\sigma^{12}]=[\sigma^{21}]\in H^2\(G, C(U_{ij},\TT)\)$,
which in particular implies that
$[\sigma^1(x)]=[\sigma^2(x)]$. Thus
$\varphi^{\alpha}$ is well defined. The continuity of
$\varphi^{\alpha}$ follows from  the continuity of the
evaluation map $x\mapsto[\sigma(x)]$ on $U$ as shown in the proof of
\lemref{lem-fix}.

Finally, suppose that $\beta=\Ad(w)\circ\alpha$ for some $1$-cocycle
$w$.   Since $\varphi^\alpha$ is defined locally, we can assume that
$\alpha=\Ad(V)$ for some $\sigma$-homomorphism $V$.  Then
$\beta=\Ad(wV)$, and it is easily checked that $wV$ is a
$\sigma$-homomorphism.
\end{proof}

Notice that if $\beta:G\to \Aut\(C_0(X,\K)\)$ is a
$C_0(X)$-linear action, then the element $\varphi^{\beta}$
constructed above is the same as the map $\varphi^{\beta}$
which appeared in Theorem~\ref{brauer}. Notice also that
$\varphi^{\alpha}=0$ if and only if all of the induced actions
$\alpha^x$ of $G$ on the fibres $A(x)$ of $A$ are unitary.

\begin{prop}\label{prop-locuni}
Let $G$ be a smooth group such that $\gab$ is compactly
generated, and let $A\in \CR(X)$. Let
$\Phi^*:C\(X,H^2(G,\TT)\)\to \E_G(X)$ denote the
splitting homomorphism for the short exact sequence
$$1\arrow{e} H^1(X,\shgab)\arrow{e}\E_G(X)\arrow{e}
C\(X,H^2(G,\TT)\)\arrow{e}
1$$
as constructed in the proof of Theorem~\ref{brauer}.
If $\alpha:G\to \Aut(A)$ is locally inner and
$[\gamma]=\Phi^*(\varphi^{\alpha})$, then
$\alpha\Xtensor{\igamma}$ is locally unitary.
\end{prop}
\begin{proof}
By the construction of $\Phi^*$ we know that there
exists a Borel cocycle $\sigma\in H^2\(G,C(X,\TT)\)$ and
a $\sigma$-homomorphism $W:G\to \UM\(C_0(X,\K)\)$ such that
${\igamma}=\Ad W$ and such that
$[\sigma(x)]=\varphi^{\alpha}(x)^{-1}$ in $H^2(G,\TT)$ for
all $x\in X$. Since $\alpha$ is locally inner, we know
further that for each $x\in X$ there exists a neighborhood
$U$ of $x$, a cocycle $\om\in Z^2\(G, C(U,\TT)\)$, and
a $\om$-homomorphism $V:G\to \UM(A_U)$ such that
$\alpha^U=\Ad V$. Then
$\varphi^{\alpha}(x)=[\om(x)]$ for all $x\in U$ by Lemma
\ref{lem-function}. It follows that the restriction
$(\alpha\Xtensor\igamma)^U$ of $\alpha\Xtensor\igamma$
to
$\(A\Xtensor C_0(X,\K)\)_U\cong
A_U\otimes_{C_0(U)}C_0(U,\K)$ is implemented by the
$\om\cdot\sigma$-homomorphism $s\mapsto V_s\otimes_UW_s$.
Since $[\om\cdot\sigma(x)]=0$ for all $x\in U$, it follows
now from \cite[Theorem 2.1]{ros2} that there exists a
neighborhood $U_1\subseteq U$ of $x$ such that
the restriction of $\om\cdot\sigma$ to $U_1$ is trivial
in $H^2\(G, C(U_1,\TT)\)$. But this implies that
the restriction of $\alpha\Xtensor{\igamma}$ to
$\(A\Xtensor C_0(X,\K)\)_{U_1}$ is unitary.
This completes the proof.
\end{proof}

\begin{thm}\label{thm-locinner}
Let $G$ be a smooth group such that $\gab$ is compactly
generated.
Suppose that $A\in \CR(X)$, and that $\alpha:G\to \Aut(A)$
is locally inner.
Then there exists a $C_0(X)$-linear action
$\beta^{\alpha}:G\to \Aut\(C_0(X,\K)\)$,
unique up to exterior equivalence, such that the stabilized
action $\alpha\Xtensor\id$ on $A\Xtensor C_0(X,\K)\ (\cong
A\otimes\K)$ is
exterior equivalent to the diagonal action $\id\Xtensor\beta^\alpha$
of $G$ on $A\Xtensor C_0(X,\K)$.

Moreover, if $\LI_G(A)$ denotes the set of exterior
equivalence classes of locally inner $G$-actions on $A$,
then $\alpha\mapsto\beta^{\alpha}$ factors through a
well defined injective map
$[\alpha]\mapsto [\beta^{\alpha}]$ of $\LI_G(A)$ into
$\E_G(X)$, which is a bijection if $A$ is stable.
\end{thm}

\begin{proof} Let $\alpha:G\to \Aut(A)$ be locally inner,
let $\Phi^*:C\(X,H^2(G,\TT)\)\to \E_G(X)$ denote
the splitting homomorphism of Theorem~\ref{brauer},
and let $[\gamma]=\Phi^*(\varphi^{\alpha})$.
Then, by  Proposition~\ref{prop-locuni},
$\alpha\Xtensor{\igamma}$ is a locally unitary action of
$G$ on $A\otimes\K$.

Let $\zeta(\alpha\Xtensor{\igamma})\in H^1(X,\shgab)$
denote the Phillips-Raeburn obstruction (see Proposition
\ref{prop-PRobs}). If $\delta:G\to \Aut\(C_0(X,\K)\)$ is locally
unitary with
$\zeta(\delta)=\zeta(\alpha\Xtensor{\igamma})$, then
it also follows from Proposition~\ref{prop-PRobs} that
$\alpha\Xtensor{\igamma}$ is exterior equivalent to
the diagonal action $\id\Xtensor\delta$ of
$G$ on $A\Xtensor C_0(X,\K)$.
Since taking diagonal actions on balanced tensor products
preserves exterior equivalence by Lemma~\ref{lem-PR3.10}, and
since ${\igamma}\Xtensor\gamma$ is exterior equivalent
to the trivial action $\id\Xtensor\id$ (Remark~\ref{rem-help}), it
follows that for
$\beta'=\delta\Xtensor\gamma$ on $A\Xtensor\coxk\Xtensor\coxk$,
\begin{align*}
\alpha\Xtensor\id\Xtensor \id&\sim
\alpha\Xtensor({\igamma}\Xtensor\gamma)
\sim(\alpha\Xtensor{\igamma})\Xtensor\gamma
\sim (\id\Xtensor\delta)\Xtensor\gamma \\
&\sim \alpha\Xtensor\beta',
\end{align*}
where $\sim$ denotes exterior equivalence.  Since
$\coxk\Xtensor\coxk\cong \coxk$, it follows that
$\alpha\Xtensor\id\sim \id\Xtensor\beta$ for some $\beta$ on $\coxk$.

We have to show that $\beta$ is unique up to exterior
equivalence. For this observe that
$\alpha\Xtensor\id\sim\id\Xtensor\beta$ implies that
$\varphi^{\alpha}=\varphi^{\alpha\Xtensor\id} =
\varphi^{\id\Xtensor\beta}=\varphi^{\beta}$ (\lemref{lem-function}).
Thus, if $\beta'$ is another $C_0(X)$-action of $G$ on
$C_0(X,\K)$ such that $\alpha\Xtensor\id\sim
\id\Xtensor\beta'$, then it follows that
$\varphi^{\beta'}=\varphi^{\alpha}=\varphi^{\beta}$.
Thus if $[\gamma]=\Phi^*(\varphi^{\alpha})$ is as
above, then
$$\id\Xtensor(\beta'\Xtensor{\igamma})\sim
\alpha\Xtensor\id\Xtensor{\igamma}\sim
\id\Xtensor(\beta\Xtensor{\igamma}),$$
which implies that the Phillips-Raeburn obstructions
$\zeta(\beta'\Xtensor{\igamma})$,
 and
$\zeta(\beta\Xtensor{\igamma})$ coincide (\lemref{lem-PR3.10}).
But then it follows from \propref{prop-PRobs} that
$\beta'\Xtensor{\igamma}\sim \beta\Xtensor{\igamma}$
which, via multiplication with $\gamma$, implies that
$\beta\sim \beta'$.

It   follows that there is a well defined
map $[\alpha]\mapsto[\beta^{\alpha}]$ from
$\LI_G(A)$ into $\E_G(A)$ which is determined by
the property that
$\alpha\Xtensor\id\sim \id\Xtensor\beta^{\alpha}$.
Since $\beta\sim\beta'$ implies that $\id\Xtensor\beta\sim
\id\Xtensor\beta'$, this map is injective.
Finally, if $A$ is stable, then we can  define
an inverse by choosing a fixed \cox-isomorphism
$\Theta: A\otimes\K\to A$ (\cite[Lemma~4.3]{pr2}), and defining
$[\beta]\to [\Ad\Theta\circ(\id\Xtensor\beta)]$
of $\E_G(X)$ onto $\LI_G(A)$.
\end{proof}

We immediately get the following corollary.

\begin{cor}\label{cor-locinner}
Let $G$ be a smooth group such that $\gab$ is
compactly generated. Let $X$ be a second countable
locally compact space and let $A\in \CR(X)$.
Then $\alpha:G\to\Aut(A)$ is locally inner if and only if
there exists $[\beta]\in \E_G(X)$ such that
the stabilized action $\alpha\Xtensor\id$ is exterior equivalent
to $\id\Xtensor\beta$.
\end{cor}

%% file: ros.tex
%
%

\section{Rosenberg's Theorem}\label{sec-appendix}

One of the important ingredients for the
proof of our results was Rosenberg's theorem
(see Theorem~\ref{thm-ros}) which implies that if
(a)~$\gab$ is compactly generated and if (b)~$H^2(G,\TT)$ is
Hausdorff, then any pointwise unitary action
on a separable continuous-trace $\cs$-algebra
$A$ is automatically locally unitary.
Our interest in smooth groups is partially explained by the fact that
all smooth groups with
$\gab$ compactly generated satisfy these
assumptions.
We give examples below which
show that neither of conditions (a)~and (b) can be weakened in
general.  On the other hand, if we assume that $A$ has continuous trace
with locally connected spectrum, then the
class of
groups with the property that pointwise unitary actions on $A$
are automatically locally unitary is significantly larger than the
class of groups which satisfy the conditions of Rosenberg's
theorem (\thmref{append}).

\begin{example}
Suppose that $G$ is a second countable locally compact abelian
group acting freely and properly on
a separable locally compact space $X$ such that $X$ is {\em not\/} a
locally trivial principal $G$-bundle.
Although Palais's Slice Theorem \cite[Theorem~4.1]{palais} implies that
$G$ cannot be a Lie group, we can, for example,
take $G=\prod_{n=1}^{\infty}\{1,-1\}$,
$X=\prod_{n=1}^{\infty}\TT$, and let $G$ act by translation on $X$.
Let $\alpha$ denote the corresponding
action of $G$ on $C_0(X)$ and let
$A=C_0(X)\rtimes_{\alpha} G$.
Then $A$ has continuous trace by \cite[Theorem 17]{green2}, and
the dual action $\widehat{\alpha} $ of $\widehat{G}$
is pointwise unitary  \cite[Proof of Theorem~3.1]{doir}.
If $\widehat{\alpha}$ would be locally unitary, then
$\spec{(C_0(X)\rtimes_{\alpha}G)\rtimes_{\widehat{\alpha}}
\widehat{G}}\cong X$
would be a locally trivial principal $G$-bundle with respect to the
double dual action
$\widehat{\widehat{\alpha}}$ \cite{pr2}.
But by the Takesaki-Takai duality theorem,
this implies that $X$ is a locally trivial
principal $G$ bundle with respect to the
original action; this contradicts our original assumption.
\end{example}

Since $\widehat G$ is discrete, $G$ is smooth and therefore
$H^2(G,\T)$ is Hausdorff.  Of course, $\widehat G$ is not compactly
generated as
required in Rosenberg's theorem.

We are now going to construct in Example~\ref{ex-3} a
pointwise unitary action
of a compactly generated group $G$ on a continuous-trace algebra $A$
which is not locally unitary.  Of course,  $H^2(G,\TT)$ will fail to be
Hausdorff.
Before we start, recall that if $N$ is a closed central
subgroup of a second countable locally compact group $G$, then the
inflation-restriction sequence
\[
H^1(N,\TT)^G \arrow{e,t}{\tg} H^2(G/N,\TT) \arrow{e,t}{\inf} H^2(G,\TT)
\]
is exact at $H^2(G/N,\TT)$, where $H^1(N,\TT)^G$ denotes the group
of $G$-invariant characters of $N$ and $\inf$ denotes
the inflation map
\cite[p.\ 53]{moore1}.

\begin{example}[cf., {\cite[p.\ 85]{moore2}}]
Our group $G$ will be a central extension of $\T^2$ by $\R^2$.
For each $\lambda \in \RR$ let $\om_{\lambda}$ denote
the two-cocycle
\[
\om_{\lambda}\((s_1,t_1),(s_2,t_2)\)=e^{i\lambda s_1t_2}
\]
on $\RR^2$. Since the real Heisenberg group is a representation group
for $\R^2$ (Example~\ref{ex-rn}), $\lambda
\mapsto [\om_{\lambda}]$ is an isomorphism between
$\RR$ and $H^2(\RR^2,\TT)$. Let
$\theta$ be any irrational number and let $\om_1$ and
$\om_{\theta}$ denote the cocycles in $Z^2(\RR^2,\TT)$ corresponding to $1$ and
$\theta$, respectively.
Let $G_1=\RR^2\times_{\om_1}\TT$ be the central extension of $\RR^2$ by
$\TT$ corresponding to $\om_1$, and let
$G=(\RR^2\times_{\om_1}\TT)\times_{\om_{\theta}}\TT$ denote the
central extension
of $G_1$ corresponding to the inflation of $\om_{\theta}$ to $G_1$.
Then $G$ is a central extension of $\RR^2$ by $\TT^2$, and is
therefore a
connected two-step nilpotent group of dimension four.
Since the cocycles involved are continuous,
$G$ is homeomorphic to
the direct product $\RR^2\times\TT^2$ with multiplication
given by
$$(s_1,t_1, z_1, w_1)(s_2,t_2,z_2,w_2)=
(s_1+s_2, t_1+t_2, e^{is_1t_2}z_1z_2,e^{i\theta s_1t_2}w_1w_2).$$
There is a natural continuous section $c$ from
$\R^2\cong G/\T^2$ onto
$G$ given by
$c(s_1,s_2):=(s_1,s_2,1,1)$.
Using the formula for the transgression map
\[
\tg:\ZZ^2\cong H^1(\TT^2,\TT)\to H^2(\RR^2,\TT),
\]
a straightforward computation shows that
$\tg(l,m)=[\om_{l+\theta m}]$.
Since $\ZZ+\theta \ZZ$ is dense in $\RR$ and $\inf$
is continuous, the identity is not closed in $H^2(G,\T)$; in other
words, $H^2(G,\T)$ is not Hausdorff.
\end{example}

\begin{example}
\label{ex-3}
We shall construct a pointwise unitary action of
the group
$G$ from the previous example
which is not locally unitary.
Let $X=\{\frac{1}{n}:n\in\NN\}\cup \{0\}$.  We define
$\alpha:G\to\Aut\(C(X,\K)\)$ as follows.
Since $\ZZ+\theta\ZZ$ is dense in
$\RR$, we find a sequence $(\lambda_n)_{n\in\NN}\subseteq
\ZZ+\theta\ZZ$
such that $\lambda_n\to 0$ in $\RR$ while
$\lambda_n\neq\lambda_m\neq 0$ for all $n,m\in \NN$, with $n\neq m$.
 Putting $\lambda_0=0$ and $\lambda_{1/n}:=\lambda_n$ we obtain
a continuous map  $x\mapsto\lambda_x$ of $X$ to $ \ZZ+\theta\ZZ\subseteq \RR$.
For each $x\in X$ let $\dot V_x:\RR_2\to \U(L^2(\RR^2)$ denote
the regular $\om_{\lambda_x}$-representation of $\RR^2$, which is given by
the formula
\[
(\dot V_{x}(s,t)\xi)(s',t')=e^{i\lambda_x s(t'-t)}\xi(s'-s,t'-t),
\]
and let $V_x:G\to \U\(L^2(\RR^2)\)$ denote the inflation of
$\dot V_x$ to $G$.
Since $x\mapsto\lambda_x$ is continuous, it follows that we obtain
a strongly continuous action $\alpha:G\to\Aut\(C(X,\K )\)$ with
$\K=\K\(L^2(\RR^2)\)$ given by defining
\[
\alpha_g(a)(x)=V_x(g)a(x)V_x(g)^*.
\]
Since $V_x$ is an $\inf(\om_{\lambda_x})$-representation for each $x\in X$
and $\alpha$ is implemented pointwise by the representations $V_x$,
and since each $[\om_{\lambda_x}]$ lies in the
range of the transgression map, it follows that $\alpha$ is pointwise unitary.

We claim that $\alpha$ is not locally unitary.
Since $X$ has only one accumulation point,
$\alpha$ is locally unitary if and only if it is unitary.
So assume that there were a strictly continuous
homomorphism $U:G\to \U\(C(X,\K)\)$ which implements $\alpha$.
Thus, for each $x\in X$ and $g\in G$ we would obtain
$$U_x(g)a(x)U_x( g)^*= V_x(g)a(x)V_x(g)^*,$$
from which it follows that
$$V_x^*(g)U_x(g)=\gamma_x(g)1$$
for some $\gamma_x(g)\in \TT$.
Since,
by construction, the maps $(x,g)\to V_x(g)$ and $ (x,g)\to U_x(g)$ are
strongly continuous,
$(x,g)\to \gamma_x(g)$ defines a continuous map $\gamma:X\times G\to\TT$.
 Moreover, since $V_x|_{\TT^2}\equiv 1$, it follows
that $\chi_x=\gamma_x|_{\TT^2}$ is a character of $\TT^2$ for all
$x\in X$.
By continuity we have $\chi_{\frac{1}{n}}\to \chi_{0}$ in
${\widehat{\T}}^2$.
Moreover, since $V_{0}$ is a unitary representation,
it follows that $V_{0}$ and $U_0$ are both unitary
representations which implement $\alpha$ at the point $0$.
But this implies that $\gamma_{0}$ is a character of $G$.
Thus, multiplying each $U_x$ with $\overline{\gamma}_{0}$, we may assume that
$U_0=V_0$. In particular, this implies
that $\chi_{0}$ is the trivial character of $\TT^2$.

We finally show that $\chi_{\frac{1}{n}}$ is not trivial for all
$n\in \NN$. Since ${\widehat{\T}}^2\cong\Z^2$ is discrete,
this will contradict
the
fact that $\chi_{\frac{1}{n}}\to\chi_{0}$.
Assume that $\chi_{\frac{1}{n}}$ is trivial for some $n\in \NN$.
Then $U_{\frac{1}{n}}|_{\TT^2}\equiv 1$, from which it follows that
$U_{\frac{1}{n}}$ is actually inflated from some unitary representation
$\dot U_{\frac{1}{n}}:\RR^2 \to\U\(L^2(\RR^2)\)$.
Since, by construction, $V_{\frac{1}{n}}$ is inflated from the regular
$\om_{\lambda_n}$-representation, say $\dot V_{\frac{1}{n}}$
of $\RR^2$, it follows that $\dot U_{\frac{1}{n}}$ and $\dot V_{\frac{1}{n}}$
implement the same action of $\RR^2$ on $\K\(L^2(\RR^2)\)$,
which contradicts the fact that $[\om_{\lambda_n}]$,
the Mackey-obstruction for the action implemented by
$V_{\frac{1}{n}}$, is non-trivial
in $H^2(\RR^2,\TT)$.
\end{example}

Note that the space $X$ in the above example is totally disconnected;
in particular, the point $0$ has no connected \nbhd s in $X$.
The following theorem shows that this lack of connectedness plays
a crucial r\^ole in our counterexample.

\begin{thm}
\label{append}
Suppose that $N$ is a closed
normal subgroup of a second countable locally compact group $G$ such that
\begin{enumerate}
\item
$G/N$ is compactly generated,
\item  $H^2(G/N,\TT)$ is Hausdorff, and
\item $H^2(N,\TT)$ is Haudorff and $\nab:=
N/\overline{[N,N]}$ is compact.
\end{enumerate}
Suppose further that $A$ is a separable continuous-trace $\cs$-algebra
such that $\hA$ is locally connected.
Then any pointwise unitary action of $G$ on $A$ is
automatically locally unitary.
\end{thm}

\begin{proof}
Let $\alpha$ be a pointwise unitary action of $G$ on $A$.
Since the properties of being
unitary and locally unitary are preserved under Morita equivalence of
systems (\cite[Proposition~3]{ech5}),
we can replace $(A,\,G,\,\alpha)$ with $(A\tensor\K,\,G,\,\alpha\tensor\id)$
and assume
that $A$ is stable.
Clearly, $(A,\,N,\,\alpha\restr N)$ is pointwise unitary, so by
Rosenberg's theorem, it is locally unitary.  Since $G$ must act trivially on
$\hA$, we can replace $A$ by an ideal and assume that $A=C_0(X,\K)$
and that $\alpha\restr N=\Ad(u)$ for a strictly continuous
homomorphism $u:N\to\UM\(C_0(X,\K)\)$.

Localizing further if necessary,
we claim that
$u$ is a Green twisting map; that is,
 $\alpha_s(u_n)=u_{sns^{-1}}$ for all $s\in G$
and $n\in N$.
Since $M\(C_0(X,\K)\)$ can be
identified with the bounded strictly continuous functions from $X$ to
$B(H)$ \cite[Corollary~3.4]{apt}, and since the strict
topology on $\UB(\H)$ coincides with the strong topology, it is not
hard to see that we may view $u$ as a strongly continuous function
from $X\times N$ to $\UB(\H)$ such that for all $n\in N$ and $x\in X$,
we have $\alpha_n(a)(x)=u(x,n)a(x)u(x,n)^*$.

In order to show that $u$ defines a Green twisting map for
$\alpha$, we need to show that
$\alpha_s\(u(\cdot,n)\)=u(\cdot,sns^{-1})$ for all $s\in G$ and $n\in
N$.
However by assumption, for each $x\in X$ there is a unitary
representation $V_x:G\to B(\H)$ such that $\alpha_s(a)(x)=V_x(s)a(x)
V_x(s)^*$ for all $a\in A$ and $s\in G$.  Since both $V_x$ and
$u(x,\cdot)$ implement the same automorphism of $\K$, there is a
character $\gamma_x$ of $\nab$ such that $u(x,n)=\gamma_x(n)V_x(n)$ for
all $n\in N$.  Now if we abuse notation slightly and write
$\alpha_s(u)(x,n)$ for $\alpha_s\(u(\cdot,n)\)(x)$, then
\begin{align*}
\alpha_s(u)(x,n) & = V_x(s)u(x,n)V_x(s)^* =
\gamma_x(n)V_x(s)V_x(n)V_x(s^{-1}) \\
&=\gamma_x(n) V_x(sns^{-1}) = \gamma_x(n)\overline{\gamma_x(sns^{-1})}
u(x,sns^{-1}).
\end{align*}
For each $x\in X$, $sN\in G/N$, and $n\in N$,
define $\lambda(x,sN)(n) = \gamma_x(n)\overline{\gamma_x(sns^{-1})} =
\gamma_x(n) \overline{s\cdot\gamma_x(n)}$,
where $s\cdot\gamma:=\gamma(s\cdot
s^{-1})$.
Clearly, $u$ will be a twisting map exactly when we can arrange for
$\lambda$ to be identically one\footnote{This invariant was also
studied \cite[\S5]{90a}.  The definition was slightly different there
--- partly to make equations such as \eqref{eq-lamb'} more
attractive
than we require here.}.
Since $A=C_0(X,\K)$, it
is not hard to see that the map
\[
(x,sN,n)\mapsto u(x,sns^{-1})\alpha_s(u)(x,n) = \lambda(x,sN)(n)1_A
\]
is continuous.  Consequently, we can view $\lambda$ as a continuous
function from $X\times G/N$ into $\hnab$.  Notice that
\begin{equation}
\label{eq-lamb'}
\lambda(x,stN) =
t\cdot\lambda (x,sN)\lambda(x,tN).
\end{equation}

Fix $x_0\in X$.  Since we may pass to still another ideal of $A$, it
 will suffice to produce a \nbhd{} $U$ of $x_0$ in
$X$ such that $\lambda(x,sN)(n)=1$ for all $sN\in G$, $n\in N$, and
$x\in U$.  Of course, replacing $u(x,n)$ by
$\overline{\gamma_{x_0}(n)}u(x,n)$, we may assume that
$\lambda(x_0,sN)(n)=1$ for all $s\in G$ and $n\in N$.
In fact, since $\hnab$ is discrete, given any $t\in G$, there is a
\nbhd{} $U_{t}\times V_t\subseteq X\times G$ of $(x_0,t)$ such that
$\lambda(x,sN)=1$ provided $(x,sN)\in U_t\times V_t$.  In view of
\eqref{eq-lamb'}, $\lambda(x,sN)=1$ for all $x\in U_t$ and $sN$ in the
subgroup of $G/N$ generated by $V_t$.  By condition~(1.), there is
 a compact set $K$ which generates $G/N$.  We can choose
$t_1,\dots,t_n$ and \nbhd s $(U_{t_1}\times V_{t_1}),\dots,(U_{t_n}\times
V_{t_n})$
such that $K\subseteq \bigcup_i V_{t_i}$.  Then we can let
$U=\bigcap_{i=1}^n U_{t_i}$.  Then $\lambda(x,sN)(n)=1$ for all $x\in U$,
$s\in G$, and $n\in N$.  Then after passing to the ideal of $A$
corresponding to $U$, we can indeed assume that $u$ is a Green
twisting map.

Now let $\rho_{x_0}$ be the element of $\hA$
corresponding to the point $x_0$ as chosen above.
Then by the above constructions, there exists
a covariant representation $(\rho_{x_0}, V_0)$ of
$(A,\,G,\,\alpha)$ such that $V_0|_N=\rho_{x_0}\circ u$,
which just means that
$(\rho_{x_0}, V_0)$ preserves the twist $u$ in the sense of Green.
Since $A$ is stable, it follows from \cite[Corollary~1]{ech5}
that $(A,\,G,\,\alpha,\,u)$ is exterior equivalent to $(A,\,G,\,\beta,\,1)$,
for some action $\beta$ of $G$ on $A$.  Then
$\beta$ is inflated from an
action $\dot\beta$ of $G/N$
(see Remark~1 on page~176 of \cite{ech5}). We are now going to show
that $\dot\beta$ is also pointwise unitary.
Since $\hA$  is
locally connected we may localize further in order to assume that
$\hA$ is connected, and we may also assume that $\hA$
is compact.
Let $(\rho_{x_0},U_0)$ denote the representation of
$(A,\,G,\,\beta)$ corresponding to $(\rho_{x_0}, V_0)$ via the
exterior equivalence between $(A,\,G,\,\alpha,\,u)$ and
$(A,\,G,\,\beta,\,1)$. Then $(\rho_{x_0},U_0)$ preserves $1$,
since $(\rho_{x_0}, V_0)$ preserves $u$. Thus it follows that
$U_0$ is inflated from a representation $\dot U_0$ of $G/N$.
Now, by assumption, $\beta$ is pointwise unitary, which implies that
$\dot\beta$ induces the trivial action of $G/N$ on $\hA$.
For each $\rho\in \hA$ let $[\om_{\rho}]\in H^2(G/N,\TT)$
denote the Mackey obstruction to extend $\rho$ to a covariant
representation of $(A,\,G/N,\,\dot\beta)$. Then $[\om_{\rho_{x_0}}]=0$
and
the map $\rho\mapsto [\om_\rho]$ is continuous by \cite[Lemma~3.3]{doir}
(or Lemma 5.3 above). Since $\hA$ is connected, it follows that
its image, say $M$, is a compact and connected subset of $H^2(G/N,\TT)$.
But since $\beta$ is pointwise unitary, it follows that
$M$ lies in the kernel of the inflation map $\inf:H^2(G/N,\TT)\to
H^2(G, \TT)$, and hence in the image of the transgression map
$\tg: H^1(N,\TT)^G \to H^2(G/N,\TT)$.
By assumption, $\nab$ is compact, so
$H^1(N,\TT)^G$
is discrete and countable (by the separability assumptions).
Thus $M$ is a  countable and connected compact Hausdorff space,
which implies that
$M$ consists of a single point.
(For example, Baire's Theorem implies that a countable compact
Hausdorff space has a clopen point.)
But since $[\om_{\rho_{x_0}}]$
is trivial, it follows that $[\om_{\rho}]$ is trivial for all
$\rho\in \hA$; in other words, $\dot\beta$ is pointwise unitary.

Now we can apply Rosenberg's theorem to the system $(A,\,G/N,\,\dot\beta)$,
from which follows that $\dot\beta$ is locally unitary.
But this implies that $\beta$, and hence also $\alpha$ is
locally unitary.
\end{proof}

We now recall that $G$ is a \fdbar{}~\emph{group} if
$\overline{[G,G]}$ is compact and $G/\overline{[G,G]}$ is abelian.
These groups are of particular interest since every
type~I \fdbar{} group has a continuous-trace group \cs-algebra
\cite[Lemma 6]{echkan} and there
are no known examples of groups with continuous-trace \cs-algebra which
are not \fdbar{} groups.

\begin{cor}\label{cor-FD}
Suppose that $G$ is a separable compactly generated
$[FD]\bar{}$-group, or that $G$ is a connected nilpotent Lie group.
Then every pointwise unitary action of $G$ on a separable
continuous-trace algebra with locally connected spectrum is
locally unitary.
\end{cor}
\begin{proof}
If $G$ is a separable compactly generated \fdbar-group,
then the theorem applies with the normal subgroup $N=\overline{[G,G]}$.
So assume that $G$ is a connected nilpotent Lie group.
Then there exists a maximal torus $T$ in the center of $G$ such that
$G/T$ is simply connected\footnote{Any connected Lie group
is a quotient of a simply connected Lie group by some central discrete group,
thus the center of a connected nilpotent group is of the form $\RR^l\times
\TT^m$ and the quotient of $G$ by $\TT^m$ is a simply connected nilpotent
group},
and
hence
$H^2(G/T,\TT)$ is Hausdorff, since
$G/T$ is smooth.
\end{proof}